\documentclass{article}

\usepackage{algorithm}
\usepackage[noend]{algpseudocode}
\usepackage{algorithmicx}

\algblockdefx{Input}{EndInput}{\textbf{Input: }}{}
\algblockdefx{Output}{EndOutput}{\textbf{Output: }}{}
\algnotext{EndInput}
\algnotext{EndOutput}

\usepackage{arxiv}

\usepackage{cite}

\newcommand{\citep}[1]{\cite{#1}}
\newcommand{\citet}[1]{\cite{#1}}

\usepackage{amsmath,amssymb,amsthm}
\usepackage{dsfont,mathtools}
\usepackage{xspace}
\usepackage{epsfig}
\usepackage{color}
\usepackage{url}
\usepackage{graphicx}
\usepackage{tikz}
\usepackage{float}
\usepackage{booktabs}
\usepackage{makecell}
\usepackage{multicol}
\usepackage{multirow}
\usepackage{csquotes}
\usepackage{wasysym}
\usepackage{longtable}
\usepackage{colortbl}
\usepackage{paralist}
\usepackage{hyperref}
\usepackage{ragged2e}

\theoremstyle{definition}
\newtheorem{definition}{Definition}[section]

\renewcommand{\today}{April 14, 2022}

\begin{document}

\title{\texorpdfstring{MOLE: Digging Tunnels Through\\Multimodal Multi-Objective Landscapes}{MOLE: Digging Tunnels Through Multimodal Multi-Objective Landscapes}}

\author{
  Lennart Sch{\"a}permeier\\
  Big Data Analytics in Transportation\\
  TU Dresden \\
  Dresden, Germany \\
  \texttt{lennart.schaepermeier@tu-dresden.de} \\
   \And
  Christian Grimme \\
  Statistics and Optimization Group\\
  University of M{\"u}nster \\
  M{\"u}nster, Germany \\
  \texttt{christian.grimme@uni-muenster.de}
  \And
  Pascal Kerschke \\
  Big Data Analytics in Transportation\\
  TU Dresden \\
  Dresden, Germany \\
  \texttt{pascal.kerschke@tu-dresden.de}
}

\maketitle

\pagestyle{plain}
\thispagestyle{fancy}
\lfoot{\vspace*{-0.875cm}\rule{\columnwidth}{0.8pt}\\
\vspace*{-0.15cm}\footnotesize \begin{justify}This version has been accepted for publication at the \textit{Genetic and Evolutionary Computation Conference (GECCO ’22), July 9–13, 2022, Boston, USA}. Permission from the authors must be obtained for all other uses, in any current or future media, including re\-printing/re\-pub\-lishing this material for advertising or promotional purposes, creating new collective works, for resale or redistribution to servers or lists, or reuse of any copyrighted component of this work in other works. \href{https://doi.org/10.1145/3512290.3528793}{DOI: 10.1145/3512290.3528793}
\end{justify}}\cfoot{}

\begin{abstract}

Recent advances in the visualization of continuous multimodal multi-objective optimization (MMMOO) landscapes brought a new perspective to their search dynamics. Locally efficient (LE) sets, often considered as traps for local search, are rarely isolated in the decision space. Rather, intersections by superposing attraction basins lead to further solution sets that at least partially contain better solutions. The Multi-Objective Gradient Sliding Algorithm (MOGSA) is an algorithmic concept developed to exploit these superpositions. While it has promising performance on many MMMOO problems with linear LE sets, closer analysis of MOGSA revealed that it does not sufficiently generalize to a wider set of test problems.
Based on a detailed analysis of shortcomings of MOGSA, we propose a new algorithm, the Multi-Objective Landscape Explorer (MOLE). It is able to efficiently model and exploit LE sets in \mbox{MMMOO} problems. An implementation of MOLE is presented for the bi-objective case, and the practicality of the approach is shown in a benchmarking experiment on the Bi-Objective BBOB testbed.

\end{abstract}

\keywords{continuous optimization, multi-objective optimization, multimodality, local search, heuristics}

\section{Introduction}

Multi-objective (MO) optimization, i.e., the optimization of multiple (conflicting) objectives at the same time, is one of the most active research areas in evolutionary computation. Many well-known algorithms like NSGA-II~\cite{Deb2002nsg2}, SMS-EMOA~\cite{BNE2007smsemoa}, or MOEA/D~\cite{ZZG2012moead} have emerged as global optimizers in many variants for approximating the optimal set of trade-off solutions for given objectives. Surprisingly, for a long time, algorithm design has mainly focused on search techniques that do not consider problem landscape characteristics of MO problems. Characteristics like multimodality were mainly considered in solution preservation approaches for diversifying the global solution set~\cite{Shir2009,Deb2008Omni,Zechman2013,Tanabe2018,liu_2018,tanabe_2019,maree2019,PengI21,JavadiM21,Tanabe2019survey}. Very few approaches also archive local structures of interest~\cite{liu2019searching}. Only recently have visualization approaches uncovered interesting local structures~\cite{kerschke2017expedition,schaepermeier2020plot}, specifically in multimodal problems, that have motivated researchers to develop local search heuristics for exploiting them~\cite{GrimmeKEPDT2019SlidingToThe,grimme2019boon}. A very interesting idea is proposed in the so-called multi-objective gradient sliding algorithm (MOGSA), which combines MO gradient descent with the discovery and traversal of locally efficient sets for reaching better local solutions sets. The algorithm relies on the observation (made in~\cite{kerschke2017expedition}) that MO attraction basins of locally efficient solutions overlap each other, which allows a (quasi automatic) descent into a better attraction basin by moving along a local solution set. In this paper, we specifically investigate components of MOGSA and highlight several shortcomings of this initial (and very simple) heuristic. We extend several aspects of the original procedure and propose a refined method called Multi-Objective Landscape Explorer (MOLE). We further empirically show that this approach -- although it is still a \emph{local} search -- is almost competitive to current (global) meta heuristics in finding \emph{globally} optimal solutions on the Bi-Objective BBOB testbed \cite{tusar2016bbobbiobj}. This highlights the potential of exploiting local structures and makes MOLE promising for application inside meta-heuristics.

After introducing the necessary notation and background in Sec.~\ref{sec:background}, MOLE is detailed in Sec.~\ref{sec:mole}. Next, we evaluate and discuss the approach in Sec.~\ref{sec:exp_discuss}, and finally conclude this work in Sec.~\ref{sec:conclude}.

\renewcommand{\vspace}[2]{}

\section{Background}
\label{sec:background}

Below, the fundamentals of and a visualization method for (continuous) MOO are presented. Using the latter, the concepts of MOGSA and MOLE will be introduced and subsequently discussed.

\subsection{Multimodality in MOO}

We consider continuous MO optimization problems (MOPs) of the form $F: \mathcal X \mapsto \mathbb R^m$ with decision space $\mathcal X \subseteq \mathbb{R}^d$ and $m$-dimensional objective vector $F(x) = (f_1(x), \ldots, f_m(x))^T$ with $f_i: \mathcal{X} \mapsto \mathbb{R}, i = 1, \dots, m$, which w.l.o.g. are being minimized. Besides $\mathcal X = \mathbb R^d$, we also consider box-constrained problems $l_i \leq x_i \leq u_i$ with lower and upper bounds $\mathbf l, \mathbf u \in \mathbb R^d$.
$F$ follows a partial order in which vector $x^*$ dominates another vector $x$, written as $x^* \prec x$, if $f_i(x^*) \leq f_i(x)$ for $i = 1, \dots, m$ and $F(x^*) \neq F(x)$. %
If for all objectives the inequalities are strict, $x^*$ is said to strictly dominate $x$. If neither $x^* \prec x$ nor $x \prec x^*$, and $x^* \neq x$, they are referred to as mutually non-dominating.
The solution for a MOP $F$ is usually the set of globally non-dominated solutions (or: Pareto set) $\mathcal S = \{x^* \ | \ x^* \in \mathcal X \land \nexists x \in \mathcal X: x \prec x^*\}$. These solutions are called Pareto optimal, and the image of $\mathcal S$ under $F$ is known as the Pareto front. 

\begin{figure}[!t]
    \centering

    \includegraphics[width=.45\columnwidth]{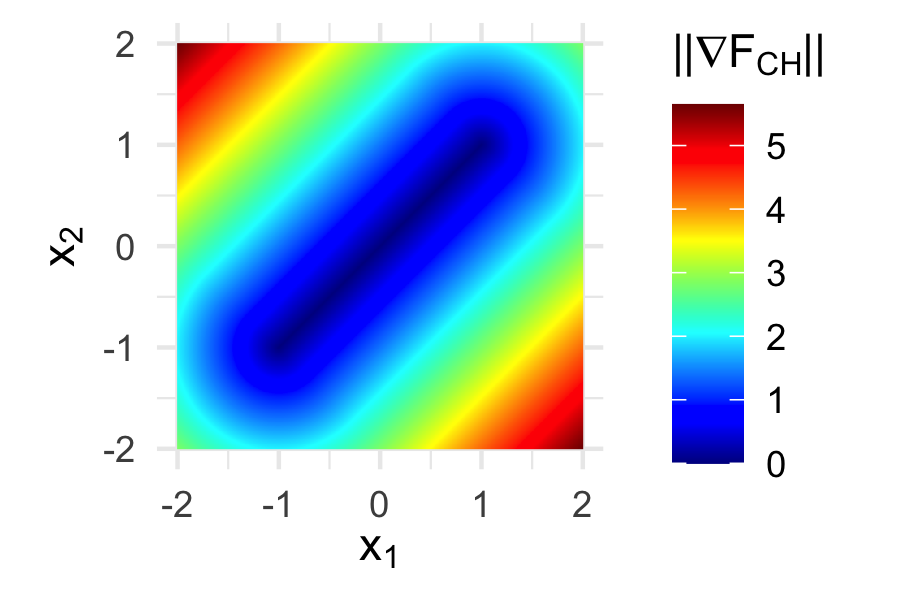}
    \includegraphics[width=.45\columnwidth]{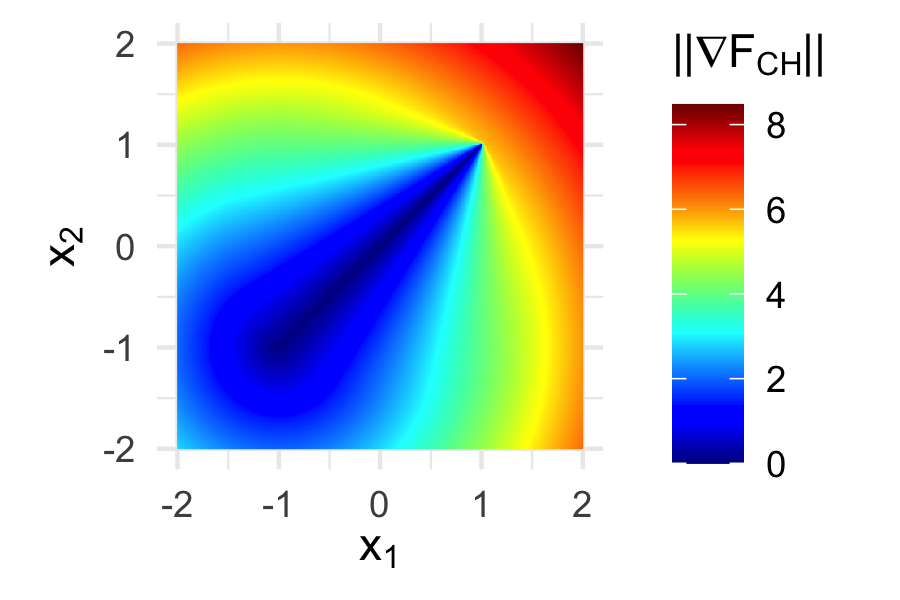}
    \vspace*{-0.55cm}
    \caption{Impact of a MOP's value ranges shown for a Bi-Sphere problem. In the right image, $f_1$ is multiplied by $100$ (compared to the left image), resulting in a strong dominance of the objective with lower range, centered at $(-1, -1)^T$.}
    \label{fig:mog-nonorm}
    
    \includegraphics[width=.45\columnwidth]{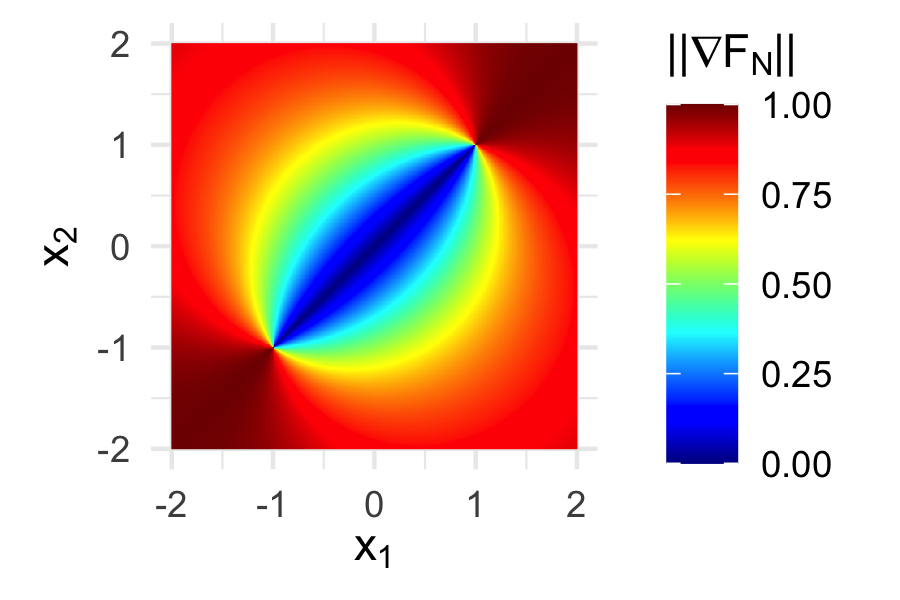}
    \includegraphics[width=.45\columnwidth]{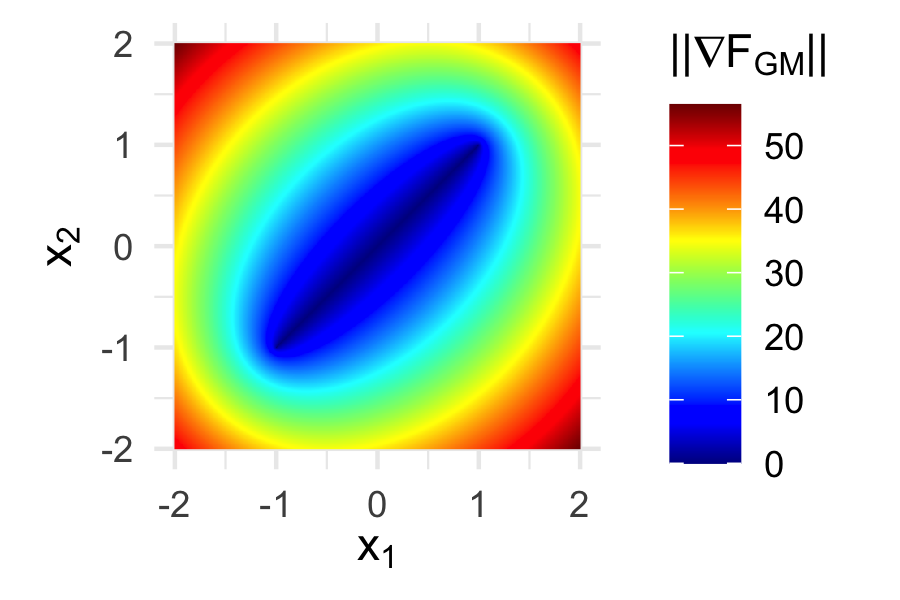}
    \vspace*{-0.55cm}
    \caption{Lengths of $\nabla F_\text{N}$ and $\nabla F_\text{GM}$, resp., on the transformed Bi-Sphere problem from Fig.~\ref{fig:mog-nonorm}. Neither landscape reveals a bias towards the objective with lower range, but in contrast to $\nabla F_\text{GM}$, $\nabla F_\text{N}$ introduces discontinuities around the SO optima.}
    \label{fig:mog-norm}
    \vspace*{-0.55cm}
\end{figure}

To discuss multimodality in continuous decision spaces, a precise notion of a point's neighborhood and local efficiency is needed. The following description is based on \citep{kerschke2016towards,grimme2019boon} and uses the Euclidean norm $\| \cdot \|$. The canonical neighborhood of a point $x \in \mathcal X \subseteq \mathbb R^d$ is the $\varepsilon$-ball $B_\varepsilon(x) = \{y \in \mathcal X \mid \|x - y\| \leq \varepsilon\}$, i.e., all points from the decision space whose Euclidean distance to $x$ is at most $\varepsilon$. Using this neighborhood, local efficiency can be defined as follows: %
\begin{definition}
    For a continuous MOP $F$, a solution $x \in \mathcal X \subseteq \mathbb R^d$ is \textbf{locally efficient (LE)} $\Leftrightarrow \exists\ \varepsilon > 0$ such that $y \nprec x$ for all $y \in B_\varepsilon(x)$.
\end{definition}
In MOPs, LE points are often connected with other LE points in their neighborhood, forming a LE set; and the set of all LE points is termed $\mathcal X_\text{LE}$. This notion is built on the topological concept of connectedness: A set $A \in \mathbb R^d$ is considered connected, if it cannot be created by the union of two disjoint (open) sets from $\mathbb R^d$. %
The connected components of a set $B$ can be defined as disjoint subsets of $B$ that are maximal w.r.t.~connectedness and whose union yields~$B$.%
\begin{definition}
    A set $\mathcal A \subseteq X$ is called a \textbf{locally efficient set}, if it is a connected component of the set of locally efficient points $\mathcal X_\text{LE}$. 
    The image of a LE set under $F$ is termed local Pareto front.
\end{definition}
There are two types of LE sets, which both affect local search dynamics differently. %
While strictly LE sets pose potential traps for local search methods, non-strictly LE sets are intersected by superposing attraction basins (see, e.g., Fig.~\ref{fig:mogsa-aspar}).
\begin{definition}
    A locally efficient set $\mathcal A$ is \textbf{strictly locally efficient}, if there is an $\varepsilon > 0$ such that no point from $\mathcal A$ is dominated by another point that is at most $\varepsilon$ distant from $\mathcal A$.
\end{definition}
\begin{figure}[t!]
    \centering
    \label{fig:strict-sets}

    \includegraphics[width=0.66\columnwidth]{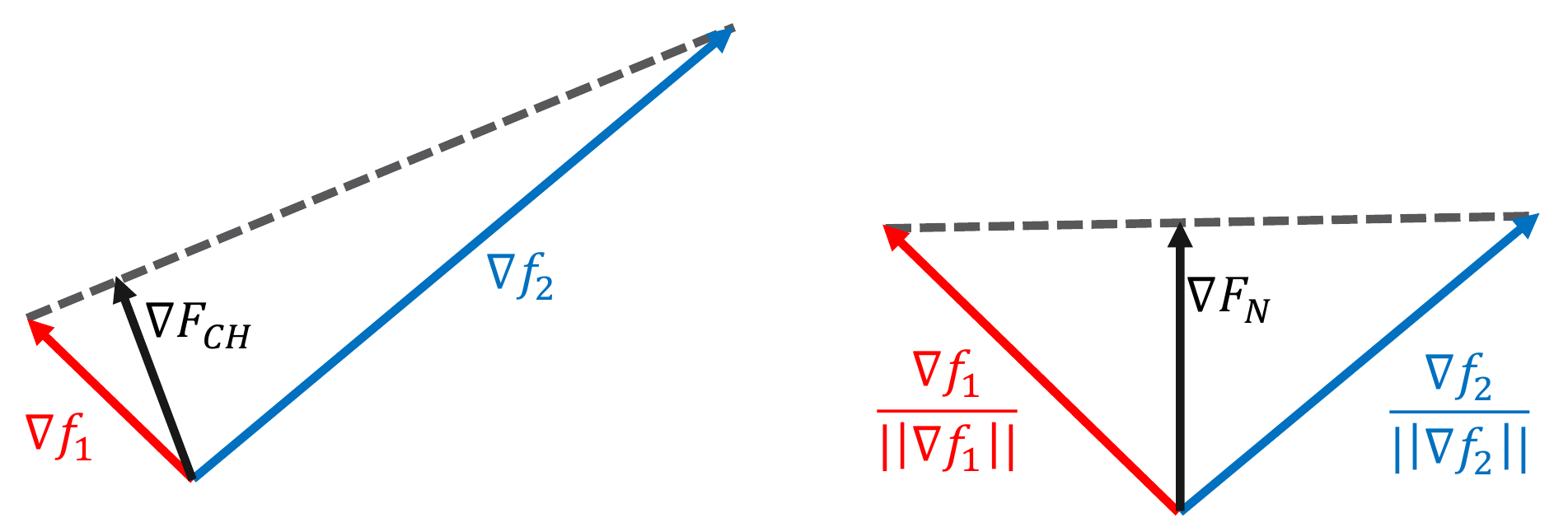}
    \vspace*{-0.2cm}
    \caption{Schematic comparison of the MOGs with and without gradient normalization. While $\nabla F_\text{CH}$ is biased towards $\nabla f_1$, $\nabla F_\text{N}$ provides an unbiased descent direction.}    
    \vspace*{-0.2cm}
    \label{fig:mog-comparison}

\end{figure}
Finally, given a continuous MOP whose individual functions $f_i$ are continuously differentiable, the first-order optimality conditions are given by their gradients:
\begin{definition}
    A decision vector $x^* \in \mathbb R^d$ of an unconstrained, continuously differentiable MOP is called first-order critical, if the Fritz-John condition holds, i.e., if there exist weights $\lambda_i \neq 0$ with $\sum_{i = 1}^k \lambda_i = 1$ such that $\sum_{i=1}^k{\lambda_i \nabla f_i(x^*)} = 0$.
\end{definition}
The concept of a multi-objective gradient (MOG) is regularly used to visualize MO landscapes and to determine the search direction and step size of MO descent algorithms in the continuous domain. Similar to gradients in the single-objective (SO) case, the (negative) MOG defines a vector field that specifies a valid descent direction for all objectives leading to a LE point where it vanishes to zero.
Definitions of the MOG, as used for MO descent, often try to generalize to the ordinary gradient when applied to the single-objective case (or to identical objectives). Such a definition based on the convex hull of the SO gradients is, e.g., given by \citet{fliege2000steepest}, as well as by \citet{desideri2012mgda}. 
Essentially, they define a MO descent direction as the shortest vector $\nabla F_\text{CH}(x) := \sum\nolimits_{i=1}^{m} \alpha^*_i \nabla f_i (x)$ within the convex hull of the SO gradients, where $\mathbf \alpha^* := \underset{\mathbf \alpha}{\mathrm{arg\,min}}\left \{\left\|\sum\nolimits_{i=1}^{m} \alpha_i \nabla f_i (x)\right\| \mid \alpha_i \geq 0, \sum\nolimits_{i=1}^{m} \alpha_i = 1 \right\}$. 
Both \cite{fliege2000steepest} and \cite{desideri2012mgda} show that this definition always leads to a valid common descent direction, if one exists. The MOG's computation is shown schematically in the left part of Fig.~\ref{fig:mog-comparison}. Unfortunately, this can also lead to an undesirable bias caused by differing lengths of the individual gradients: A SO gradient with a smaller norm has a greater influence on the selected descent direction, as shown in Fig.~\ref{fig:mog-nonorm}. Although this is not problematic for finding valid descent directions, MOPs with strongly differing values for the objectives will arbitrarily favor the one with a smaller range (i.e., shorter gradient). Thus, for a simple multiplicative transformation of a single objective, the descent direction changes, whereas the visualization of the objective space and the order of the solution sets determined by the HV indicator (cf.~Sec.~\ref{sec:mo-descent}) remain unaffected. Further, the resulting HV values would merely be affected by a scalar factor.%

As a remedy, the SO gradients are often normalized to unit vectors before applying the above MOG computation. In the bi-objective case this simplifies to the following:
\begin{equation}
    \nabla F_\text{N}(x) := 0.5 \cdot \left( \nabla f_1 (x) / \|\nabla f_1 (x)\| + \nabla f_2 (x) / \|\nabla f_2 (x)\| \right)
\end{equation}
This yields an unbiased descent direction -- independent of the lengths of the SO gradients -- but introduces a discontinuity around SO critical points, i.e., where (at least) one of the SO gradients vanishes. As shown in Fig.~\ref{fig:mog-norm} for a Bi-Sphere problem, surrounding the SO optima, there are directions which result in maximal (length 1) and minimal (length 0) values for $\nabla F_\text{N}$. This complicates the usage of $\nabla F_\text{N}$, especially if one is interested in using its length to motivate the step size during descent or if one wants to use a MO descent around the SO optima. Further, $\nabla F_\text{N}$ does not generalize the SO gradient but rather the normalized SO gradient.

$\nabla F_\text{N}$ can be modified to produce an unbiased and smooth MOG, which to the authors' best knowledge presents a novelty. Rather than normalizing each gradient in isolation, they are normalized to the geometric mean of the individual (SO) gradient lengths, i.e.:
\begin{equation}
    \nabla f_i' (x) := {\sqrt[m]{\prod\nolimits_{j=1}^m\|\nabla f_j (x)\|}} \cdot \nabla f_i (x) / \|\nabla f_i (x)\|
\end{equation}
Then, the MOG based on the transformed gradients is computed, which is equivalent to $\nabla F_\text{N}$ scaled by the geometric mean of the gradient lengths. For the bi-objective case, this yields:
\begin{equation}
    \begin{split}
        \nabla F_\text{GM} (x)
        &:= \frac 1 2 \left(
            \sqrt{\frac {\|\nabla f_2 (x)\|}{\|\nabla f_1 (x)\|}} \nabla f_1 (x) +
            \sqrt{\frac {\|\nabla f_1 (x)\|}{\|\nabla f_2 (x)\|}} \nabla f_2 (x)
        \right) \\
        &= \nabla F_\text{N}(x) \cdot \sqrt{\|\nabla f_1 (x)\|} \sqrt{\|\nabla f_2 (x)\|}
    \end{split}
\end{equation}
This enables $\nabla F_\text{GM}$ to maintain continuity around SO optima: If any $\nabla f_i(x) \rightarrow 0$, the geometric mean of the gradient lengths ensures that the MOG near this optimum also converges to zero. So, if any gradient becomes zero at some point, the MOG can be defined as zero at these points without introducing new discontinuities.
It can also be verified that it generalizes the gradient in the SO case, and that a multiplicative transformation of any objective by a factor $\gamma > 0$ results in scaling the MOG by $\sqrt[m]{\gamma}$ without affecting its direction (cf. Fig.~\ref{fig:mog-norm}).

\subsection{Visualization for Multimodal MOPs}

Visualizing continuous MOPs is helpful for understanding %
structural properties 
and for motivating %
new approaches. Here, we outline a recent visualization technique, which we use for our purposes.

The so-called \emph{Plot of the Landscape with Optimal Trade-offs (PLOT)} \cite{schaepermeier2020plot} %
enhances a previous technique, called 
gradient field heatmaps (GFH) \cite{kerschke2017expedition}. The latter visualizes the local search dynamics and LE solutions of a MOP. After discretizing the decision space into a grid, the MOG is estimated per cell. Then, a MO gradient descent is simulated which follows the path along the closest grid neighbor (w.r.t.~angle) in the respective descent direction. This process terminates, if either $\|\nabla F\| < \varepsilon$ (for small $\varepsilon > 0$), or the descent reaches a previously visited point. 
For the visualization, the lengths of all MO gradient steps along the descent path of a given point are cumulated. This value is interpreted as height value for the considered starting point. In contrast to GFH, the used PLOT approach shows this height values in gray-scale for a representation of the attraction basins of LE sets. %
The visualization of the LE sets is decoupled and depicts 
their log-scaled domination counts 
(applying a method by \citep{fonseca1995multiobjective} to $\mathcal{X}_{LE}$ only)
on a separate color scale. This technique simplifies the identification of a MOP's (close to) globally efficient regions and a visual ranking of the locality relation of LE sets. %

\subsection{MOGSA}

Based on observations in the GFH visualizations, \cite{grimme2019boon} introduced a prototypical MO local search algorithm, the \textbf{Multi-Objective Gradient Sliding Algorithm (MOGSA)}, which aims to exploit the superpositions between LE sets. Algorithm~\ref{alg:mogsa} gives a high-level overview of MOGSA's implementation \cite{kerschke2019mogsa} for bi-objective MOPs.%

\begin{algorithm}[!t]
    \caption{Multi-Objective Gradient Sliding Algorithm: MOGSA}\label{alg:mogsa}
    \begin{algorithmic}[1]
        
        \Input {MOP $F: \mathcal X \rightarrow \mathbb R^m$, starting point $x \in \mathcal X \subseteq \mathbb{R}^d$}
        \EndInput
        
        \Output{Archive of visited points}
        \EndOutput
        
        \Procedure{MOGSA}{$F, x$}
            
            \State $x_\text{next} \gets x$

            \While{$x_\text{next} \neq \texttt{NULL}$}
                \State $x^* \gets \Call{FindLocallyEfficientPoint}{x_\text{next}, F}$
                \State $x_\text{next} \gets \Call{ExploreEfficientSetMOGSA}{x^*, F}$
            \EndWhile
            
            \State \textbf{return} Archive of all visited points
    \EndProcedure
    \end{algorithmic}
\end{algorithm}
\vspace*{-0.5cm}
MOGSA proceeds as follows: Starting in a point $x \in \mathcal X$ a MO descent is applied (\textsc{FindLocallyEfficientPoint}) to reach a LE point. From the starting point, one moves downhill by $\nabla F_\text{N}$ (multiplied by a constant step size), terminating once $\nabla F_\text{N}$ %
shrinks below a (pre-defined) size, or a maximum number of iterations is reached. 
Then, \textsc{ExploreEfficientSetMOGSA} explores the discovered LE set by following the SO gradients successively for each objective with a fixed step size. If from one step to the next the angle of the tracked objective's gradient changes by over $90$°, or any of the gradients vanishes, a SO optimum is presumed. However, if the angle between the two SO gradients becomes less than $90$°, a superposed basin is presumed and the current point is returned. %
This is repeated until a strictly LE set is found, which is presumed to be globally optimal in many cases.

\section{MOLE}
\label{sec:mole}

\begin{figure}[!t]
    \centering
    
    \includegraphics[width=.425\columnwidth, trim = 0mm 0mm 0mm 62mm, clip]{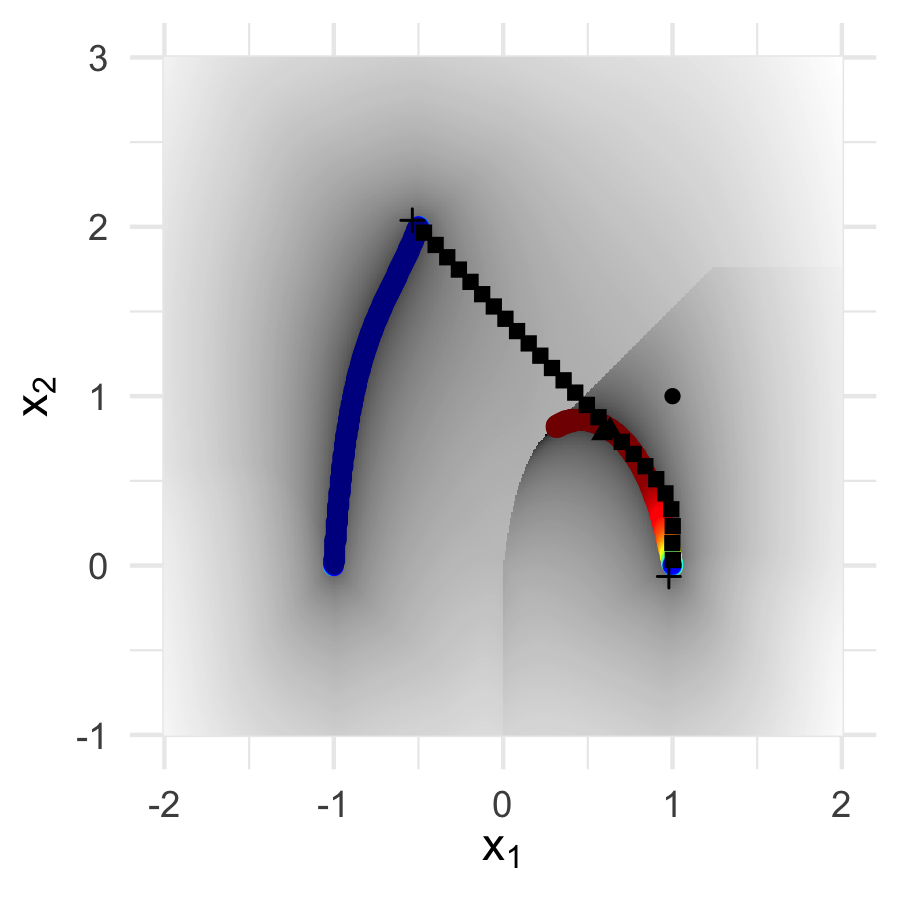}
    \vspace*{-0.2cm}
    \includegraphics[width=.425\columnwidth, trim = 0mm 0mm 0mm 62mm, clip]{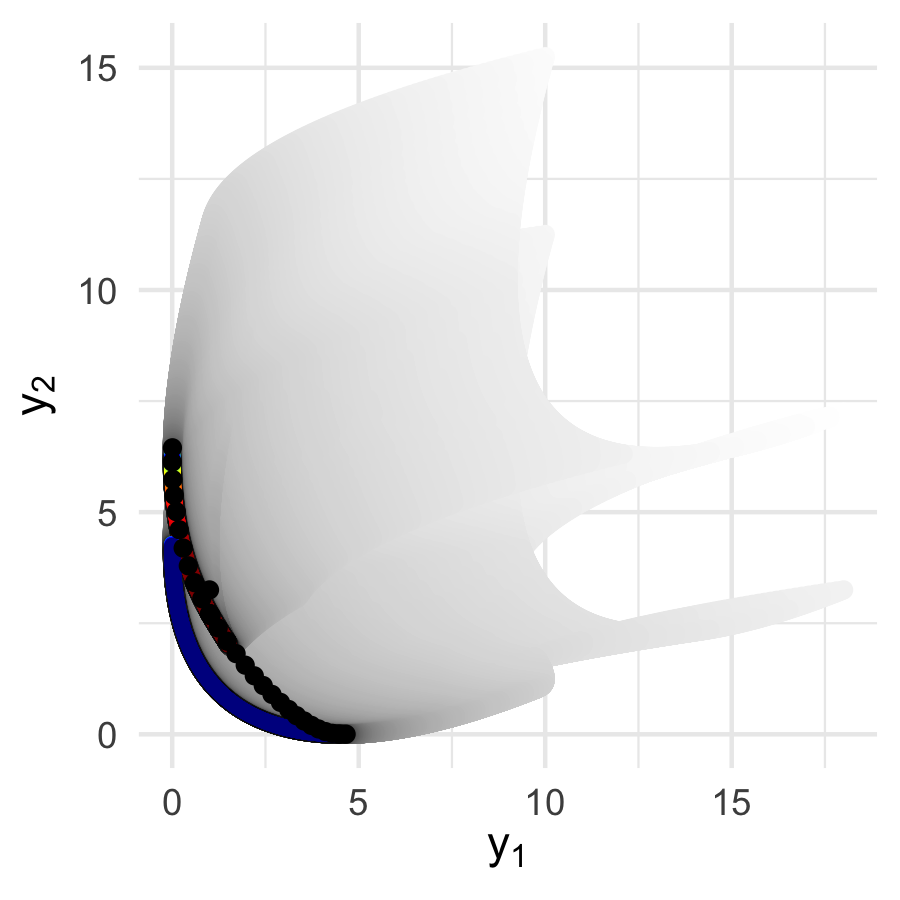}
    \vspace*{-0.3cm}
    \caption{
    Trace of a MOGSA run on the Aspar function $F(x_1, x_2) := (x_1^4 - 2 x_1^2 + 2 x_2^2 + 1, (x_1 + 0.5)^2 + (x_2 - 2)^2)$ starting in $(1,1)^T$, with the PLOT in the background (left: decision, right: objective space). %
    } %
    \label{fig:mogsa-aspar}
    
    \vspace*{-0.08cm}
    \includegraphics[width=.425\columnwidth, trim = 0mm 30mm 0mm 3mm, clip]{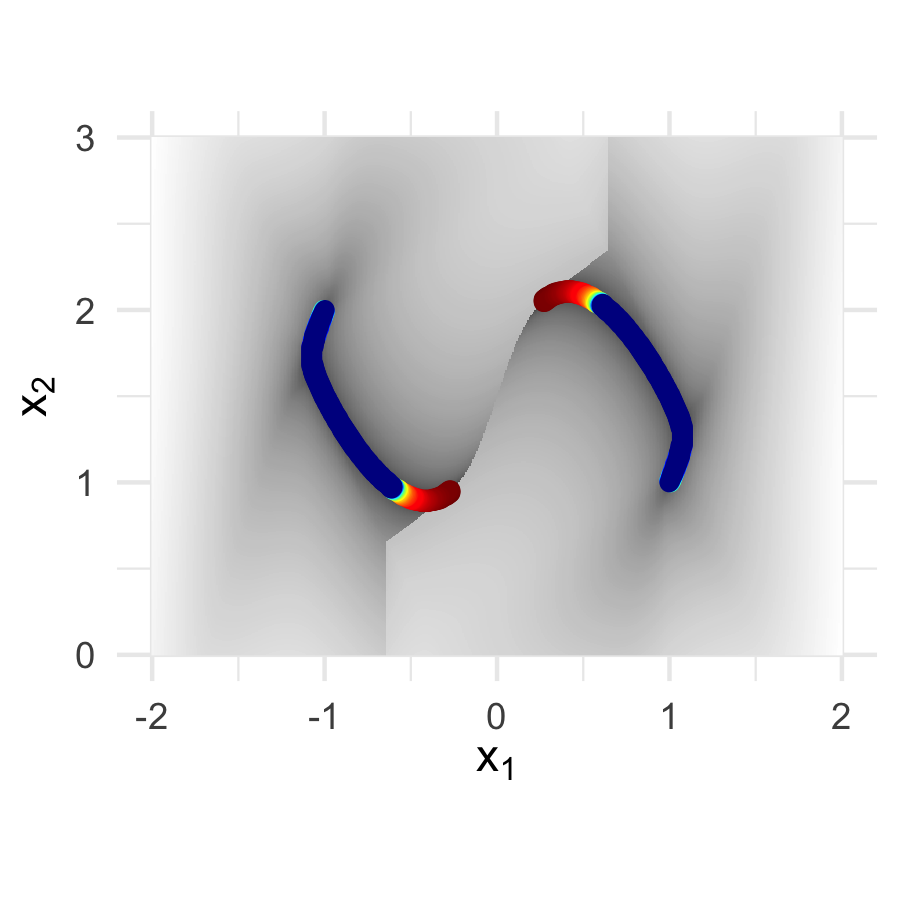}
    \includegraphics[width=.425\columnwidth, trim = 0mm 0mm 0mm 78mm, clip]{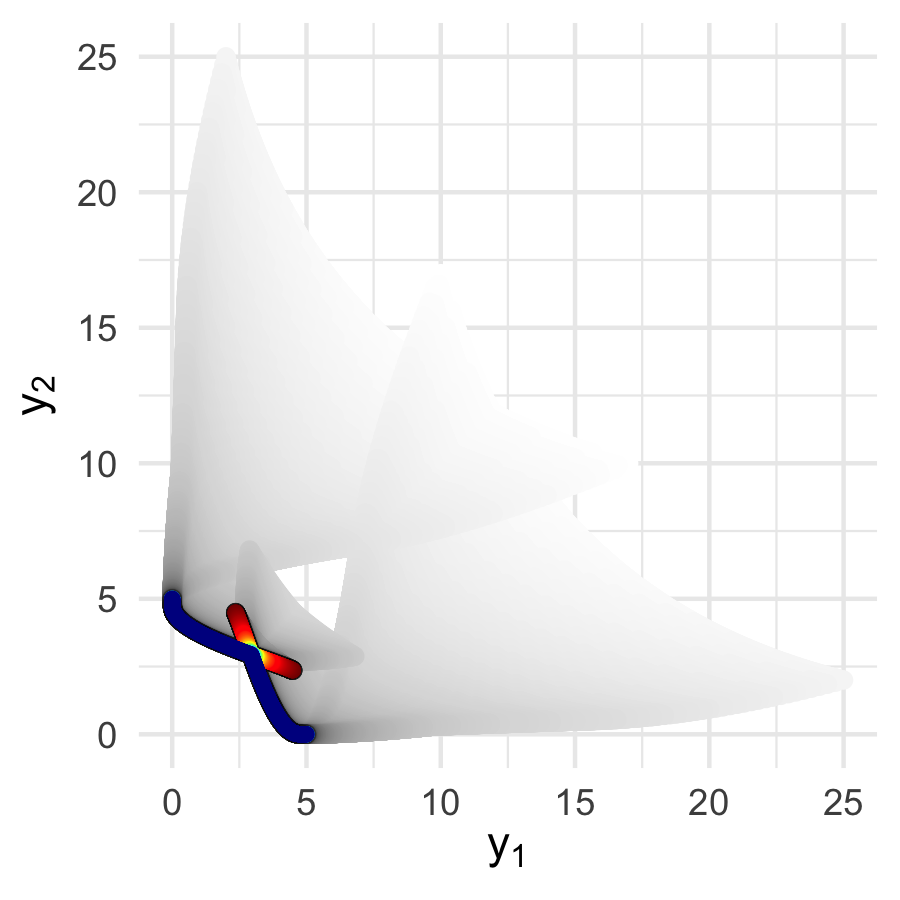}
    \vspace*{-0.47cm}
    \caption{
    PLOT of a Bi-Rosenbrock function $F(x_1,x_2) := ((1 - x_1)^2 + 1 \cdot (x_2 - x_1^2)^2, (1 + x_1)^2 + 1 \cdot (-(x_2 - 3) - x_1^2)^2)$. Instead of having a strictly LE set, it consists of two intersected sets.
    }
    \label{fig:birosen-plot}
    \vspace*{-0.35cm}
\end{figure}

Despite the great potential of the basic concept behind MOGSA, its view of the MO landscape is still quite simplified, resulting in undesirable behavior as, e.g., demonstrated in Fig.~\ref{fig:mogsa-aspar}. 
One of the weaknesses of MOGSA is its termination criterion: It stops once a strictly LE set has been found. However, the landscape of any given continuous MOP is not guaranteed to have such a set, as exemplarily shown for the Bi-Rosenbrock function in Fig.~\ref{fig:birosen-plot}. %
In fact, the globally efficient set may even span across multiple (not necessarily strict) sets, i.e., they may be (partially) dominated by a neighboring set. MOGSA would thus need to keep track of the location of the visited sets in the decision and objective space.

Once a LE point has been found, MOGSA explores the respective set by following the SO gradient descent (per objective) until either an endpoint (i.e., SO optimum) is reached, or a ridge is detected by the MOG. This presumes that the SO gradients point parallel to the LE set. However, this does not have to be the case, as demonstrated by MOGSA's behavior in Fig.~\ref{fig:mogsa-aspar}. It also demonstrates the vulnerability of the angle-based detection of a new attraction basin.

Finally, MOGSA's MO descent is rather basic. It uses $\nabla F_\text{N}$ as descent direction and its length $\|\nabla F_\text{N}\|$ determines the size of the next step. This makes the descent sensitive to the chosen step size parametrization, and thus may fail to converge onto a LE point.

\subsection{Overview}

Based on these observations, Algorithm~\ref{alg:mole} presents the  \textbf{Multi-Objective Landscape Explorer (MOLE)}, an algorithm that explores the LE sets of continuous MO optimization landscapes. Given a continuous MOP and a set of starting points $\mathcal P \subseteq \mathcal X$, MOLE 
approximates 
the LE sets and tracks the transitions between them. %
MOLE consists of several components, 
which together enable it to efficiently explore the network of LE sets of a continuous MOP. %

First, a strategy is needed for generating a suitable set of starting points $\mathcal P \subseteq \mathcal X$. These may be sampled uniformly at random, from a specific subspace, or using a different sampling strategy.

Then, a \textsc{MultiObjectiveDescent} is applied to the starting point to receive a LE point. The output should contain a point dominating its input. If no descent is possible, the input is returned. Further, the individual steps should be small enough to ensure that the starting and descended point belong to the same attraction basin.

After finding a LE point, a new search along the corresponding LE set and any further sets with superposing attraction basins is started. Here, LE points that belong to already visited sets need to be detected and handled separately to avoid looping between sets that partially dominate %
each other (cf. Fig.~\ref{fig:birosen-plot}). This is implemented in the \textsc{FindContainingSet} procedure. If no corresponding set is found, the LE set is explored using \textsc{ExploreEfficientSet}, which should return the %
approximated LE set, as well as LE points that belong to superposing attraction basins found during this process.
The superposing points are then appended to the list of points to explore during further iterations. By maintaining this list, we ensure that all LE sets, reachable from a given starting point, are explored. If there are no further points to explore, the evaluation of that starting point is finished, and the next iteration is started.

\begin{algorithm}[!t]
\small
    \caption{Multi-Objective Landscape Explorer (MOLE)}\label{alg:mole}
    \begin{algorithmic}[1]
        
        \Input {MOP $F: \mathcal X \rightarrow \mathbb R^m, \mathcal X \subseteq \mathbb R^d$, starting points $\mathcal P \subseteq \mathcal X$}
        \EndInput
        
        \Output{Approximations to the LE sets of $F$}
        \EndOutput
        
        \Procedure{MOLE}{$F, \mathcal P$}
            
            \State $\texttt{sets} \gets []$

            \While{\emph{termination criteria not met}}
                \State $p \gets \text{Next point from } \mathcal P$
                \State $p^* \gets \Call{MultiObjectiveDescent}{p, F}$
                \State $\texttt{points\_to\_explore} \gets [p^*]$
                
                \While{$\texttt{points\_to\_explore}$ \emph{is not empty}}
                    \State $x^* \gets \texttt{points\_to\_explore}.pop\_back()$
                    
                    \State $\texttt{found\_set} \gets \Call{FindContainingSet}{x^*, \texttt{sets}}$
                    
                    \If{$\texttt{found\_set} \neq \texttt{NULL}$}
                        \State Insert $x^*$ into $\texttt{found\_set}$
                    \Else
                        \State $(\texttt{set}, \texttt{super}) \gets \Call{ExploreEfficientSet}{x^*, F}$
                        
                        \State $\texttt{sets}.append(\texttt{set})$

                        \State $\texttt{points\_to\_explore}.append\_all(\texttt{super})$
                    \EndIf
                \EndWhile
                \State $\Call{PostProcess}{\texttt{sets}}$
            \EndWhile
        \State \textbf{return} $\texttt{sets}$
    \EndProcedure
    \end{algorithmic}
\end{algorithm}

Finally, a post-processing phase (\textsc{PostProcess}) is applied periodically. This component allows to adapt the MOLE algorithm to the optimization task at hand, e.g., by regularly optimizing for a quality indicator, or simply validating the LE set models. %

A run of MOLE can be terminated based on a variety of conditions depending on the application scenario, e.g., when all given starting points are evaluated or the available budget runs out.

Note that this description of MOLE is independent of the number of objectives. However, so far the implementation details are missing, and as observed with the MOGSA prototype, these can be crucial in the realization of a practical algorithm. To close this gap, we will provide necessary implementation details for MOLE in the (continuous) bi-objective case. %

\subsection{Multi-Objective Descent}\label{sec:mo-descent}

A multitude of MO descent procedures have been proposed in the literature (see \cite{fukuda2014survey} for an overview). %
Many of these methods rely on a MOG to find a descent direction that decreases all objectives simultaneously \citep{fliege2000steepest,desideri2012mgda}, and they may be extended to (Quasi-)Newton methods that exploit second-order information in the MOG vector field. %
Finally, a MO descent could also be implemented by solving an appropriate scalarized problem, e.g., by optimizing the dominated Hypervolume (HV) \cite{zitzler2003performance}, i.e., the hyperspace covered by the solutions of the approximation set (in the objective space) w.r.t.~an anti-optimal reference point, to the starting point.%

As any descent method, MOLE's MO descent requires a search direction -- it uses $d^{(t)}= -\nabla F_\text{GM}(x^{(t)})$ -- and a step size $\alpha^{(t)}$ at each time $t$ to find $x^{(t+1)}$ from the previous point $x^{(t)}$. Although not gaining the speed up that (Quasi-)Newton methods often enjoy, this approach stays consistent to the gradient-based visualizations introduced above, and is most straightforward in its implementation. Rather than selecting a more sophisticated search direction, good practical performance is ensured by the choice of step size. %

For this purpose, we essentially rely on the so-called Barzilai-Borwein (BB) step size rule~\cite{barzilai1988two}, but consider the augmented version presented by \cite{dai2014positive}. Given the last two iterates $x^{(t)}, x^{(t-1)}$ and their respective gradients $d^{(t)}, d^{(t-1)}$, the step size can be computed as
\begin{equation}\label{eq:alpha_both}
    \max \{ (s^{(t-1)T} s^{(t-1)}) / (s^{(t-1)T} y^{(t-1)}), \|s^{(t-1)}\| / \|y^{(t-1)}\| \}
\end{equation}
with $s^{(t-1)} := x^{(t)} - x^{(t-1)}$ and $y^{(t-1)} := d^{(t)} - d^{(t-1)}$.
This approach is easily extended to the MO domain: \cite{barzilai1988two} already mention that their method is not limited to using gradients. Thus, we can apply the BB step size using the MOG as the search direction without further changes. The BB step sizes were also used by \citep{morovati2016barzilai} for MO descent, although their approach uses a different search direction.
Alas, while the BB step sizes present good overall performance for very little effort (computationally and in terms of function evaluations), they only guarantee convergence in special situations. %

In SO optimization, step sizes that do not lead to oscillations around local optima can, for instance, be verified by the Armijo rule  \citep{armijo1966minimization}. Formally, it is defined with a constant $0 < \beta < 1$ (usually a small value like $10^{-4}$) and a minimization problem $f: \mathcal X \rightarrow \mathbb R$ as: %
\begin{equation}
    f(x + \alpha \cdot d) \leq f(x) + \beta \cdot \alpha \cdot d^\texttt{T} \nabla f(x)
\end{equation}
That is, a step size of $\alpha$ is accepted, if moving into the direction $d$ ensures that the function value is at least decreased by a factor of $\beta$ times the improvement that would be expected by the local linearization using the gradient. %
This concept can directly be extended to MOO by ensuring that the Armijo condition holds for each objective individually while choosing the MOG as the descent direction. %
An Armijo-like rule for the MO case is also used in \cite{fliege2000steepest}.

Verifying the Armijo condition in itself does not yet improve the performance of the BB method, as it benefits from sometimes taking larger steps, even if the objective value(s) worsen. In SO optimization, a common way to counter this problem is a %
nonmonotone line search, such as the global BB method \cite{raydan1997barzilai}.
Instead of ensuring that the Armijo condition holds when compared to the last iterate, the global BB compares it to the respective maximal value of the last few iterations, which allows some steps to worsen when compared to their predecessor, but still ensures convergence over time. %

\begin{algorithm}[!t]
\small
    \caption{Nonmonotone MO Descent with BB Stepsize}\label{alg:mo_descent}
    \begin{algorithmic}[1]
        \Input{Starting point $x \in \mathcal X$, continuous MOP $F: \mathcal X \rightarrow \mathbb R^m$}
        \EndInput
        
        \Output{LE point $x^* \in \mathcal X$}
        \EndOutput
        \Procedure{MultiObjectiveDescent}{$x, F$}
            
            \State $x^{(0)} \gets x$, $x^{(1)} \gets x$, $\alpha^{(0)} \gets \alpha_\text{min} / {\|\nabla F^{(0)}\|}$
            
            \While {($F(x^{(0)} - \alpha^{(0)} \cdot \nabla F^{(0)}) \preceq F^{(1)}$ \textbf{and} \\ \hspace{1.3cm} $\alpha^{(0)} \leq \alpha_\text{max} / \|\nabla F^{(0)}\|$)}
                \State $x^{(1)} \gets x^{(0)} - \alpha^{(0)} \cdot \nabla F^{(0)}$
                \State $\alpha^{(0)} \gets \lambda \cdot \alpha^{(0)}$
            \EndWhile
            
            \If{$x^{(0)} = x^{(1)}$}
                \State \textbf{break}
            \EndIf
            
            \For{$t=1,2,\dots,\texttt{max\_iter}$}
                \If{$\|\nabla F^{(t)}\| < \gamma$}
                    \State \textbf{break}
                \EndIf
                \State $\forall i: f^\text{ref}_i \gets \max\{f_i^{(\max(0, t - H+1))}, \dots, f_i^{(t)}\}$

                \State $\alpha^{(t)} \gets \|s^{(t-1)}\| \; / \; \min \{s^{(t-1)T} y^{(t-1)}, \; \|y^{(t-1)}\|\}$ see (\ref{eq:alpha_both})
                \State Ensure $\alpha_\text{min} / {\|\nabla F^{(t)}\|} \leq \alpha^{(t)} \leq \alpha_\text{max} / {\|\nabla F^{(t)}\|}$
                
                \State $x^{(t+1)} \gets x^{(t)} - \alpha^{(t)} \cdot \nabla F^{(t)}$
                
                \While{($\exists i: f_i^{(t+1)} \geq f^\text{ref}_i + \beta \cdot \alpha^{(t)} \cdot \nabla F^{(t)T} \nabla f_i^{(t)}$ \textbf{and} \\ \hspace{1.7cm} $\alpha^{(t)} > \alpha_\text{min} / \|\nabla F^{(t)}\|$)}
                    \State $\alpha^{(t)} \gets \max \{\alpha^{(t)} / \lambda, \alpha_\text{min} / \|\nabla F^{(t)}\|\}$
                    \State $x^{(t+1)} \gets x^{(t)} - \alpha^{(t)} \cdot \nabla F^{(t)}$
                \EndWhile

                \If{($\alpha^{(t)} = \alpha_\text{min} / \|\nabla F^{(t)}\|$ \textbf{and} \\ \hspace{1.2cm} ($x^{(t+1)} \nprec x^{(t)}$ \textbf{or} \\ \hspace{1.2cm} $\exists i: f_i^{(t+1)} \geq f^\text{ref}_i + \beta \cdot \alpha^{(t)} \cdot \nabla F^{(t)T} \nabla f_i^{(t)}$))}
                    \State \textbf{break}
                \EndIf
            \EndFor
            
            \State \textbf{return} $x^{(t)}$
    \EndProcedure
    \end{algorithmic}
\end{algorithm}

\begin{table}[!t]
    \caption{Parameters for the MO descent, as well as their suggested default values. \texttt{diag} refers to the diagonal of the decision space (if lower and upper bounds are known), or of a similar central part of the search space otherwise.} %
    \centering
    \vspace*{-0.2cm}
    \small
    \begin{tabular}{|l|l||r||l|}
        \hline
        Parameter & Value & Description \\
        \hline
        \hline
        $\gamma$ & \texttt{1e-6} & Min. length of MOG \\
        $\alpha_\text{min}$ & \texttt{1e-6} & Min. descent step size \\
        $\alpha_\text{max}$ & \texttt{diag/100} & Max. descent step size \\
        $\lambda_\text{descent}$ & \texttt{2} & Descent step scaling factor \\
        $\beta$ & \texttt{1e-4} & Factor for Armijo condition\\
        $H$ & \texttt{100} & History size for nonmon.~search\\
        \texttt{max\_iter} & \texttt{1000} & Max. number of descent iterations\\
        \hline
    \end{tabular}
    \label{tab:parameters-descent}
    \vspace*{-0.25cm}
\end{table}

Finally, the MO descent implemented for MOLE is summarized in Alg.~\ref{alg:mo_descent} and the suggested default parameters are listed in Tab.~\ref{tab:parameters-descent}.
Notations such as $F^{(t)}$ are abbreviations for $F(x^{(t)})$, and the gradients are estimated using (central) finite differences with $\epsilon = 10^{-8}$.

The descent begins with a doubling line search to find an initial step size, where the step size is repeatedly increased by a factor of $\lambda_\text{descent}$ until no dominating point is found anymore. To ensure that the steps taken are sufficiently local, a maximal step size $\alpha_\text{max}$ for an individual step is defined. In addition, a minimal step size $\alpha_\text{min}$ ensures that the resulting steps are not too small. $\alpha_\text{min}$ can also be used to verify convergence early: If we cannot find a dominating point with the minimal step size, we can terminate the descent.

Then, as long as $\|\nabla F^{(t)}\| \geq \gamma$, but for at most \texttt{max\_iter} iterations, we perform the global BB procedure discussed above. Note that the absolute step sizes are clamped to $[\alpha_\text{min}, \alpha_\text{max}]$. If the Armijo condition cannot be verified, $\alpha^{(t)}$ is reduced by a factor of $\lambda_\text{descent}$ successively. In addition to terminating when a step with size $\alpha_\text{min}$ does not lead to a point dominating the previous one, we terminate if the Armijo condition w.r.t.~the reference point cannot be verified for that step size. %
The history size of the nonmonotone search, $H$, is intentionally set high by default in order to interfere with the chosen step sizes only for searches that would fail to converge according to one of the other conditions.

\subsection{Multi-Objective Continuation}

Given a LE point, the next task of MOLE is to explore the corresponding LE set. %
In general, the abstract problem is called numerical continuation or numerical path following and is regularly addressed in MO optimization~\cite{hillermeier2001nonlinear, harada2007uniform, martin2018pareto}. 
A central continuation approach is the \textbf{predictor-corrector (PC)} method, which consists of two steps: First, a new point is predicted into a direction of interest along the desired set. %
Clearly, the resulting point may be significantly off the true set. This is subsequently corrected by numerically solving an optimization problem, which terminates at points of the traced set. %

\begin{figure}[!t]
    \centering
    \includegraphics[width=0.8\columnwidth]{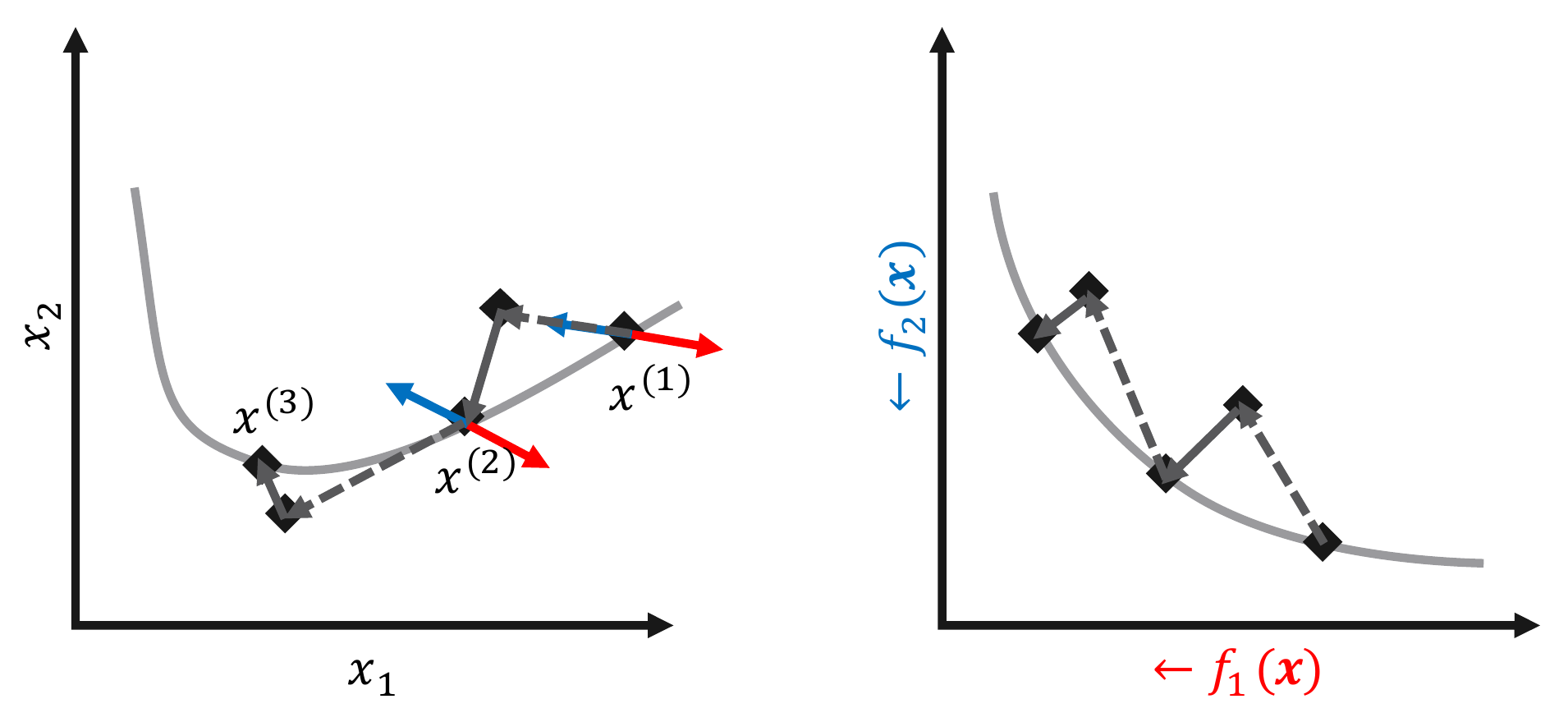}
    \vspace*{-0.25cm}
    \caption{Illustration of the proposed bi-objective continuation method in decision and objective space. For the first prediction, the negative gradient of the objective that is pursued (here: $f_1$, red) is taken. For further predictions, the secant through the last two points on the set is used.}
    \vspace*{-0.2cm}
    \label{fig:mo-continuation}
\end{figure}

In the bi-objective case, we can find a good approximation of the set direction if we have found at least two points along the set: Then, we take the secant through the last two points %
and predict that the set is going to continue in this direction \citep{martin2018pareto}. This results in a ``free'' and reasonably good prediction direction, especially if the set is (almost) linear.
As corrector, we use the previously discussed MO descent. %
MOLE's predictor and corrector are sketched in Fig.~\ref{fig:mo-continuation}. %

As MOLE is expected to work on a wide range of functions, %
an effective step size adaptation within the PC approach is needed. This also serves the purpose of efficiently exploring (almost) linear parts of the LE sets, while focusing on areas with stronger curvature. Further, PC methods can only describe a single LE set; yet, we are also interested in detecting superposing attraction basins during exploration and thus also need to detect the transitions between sets. 
The resulting procedure to explore a LE set with MOLE is outlined in Alg.~\ref{alg:explore-set} and its default parameters are given in Tab.~\ref{tab:parameters-explore-set}.

\begin{algorithm}[!t]
\small
    \caption{Exploration of a Bi-objective LE Set}\label{alg:explore-set}
    \begin{algorithmic}[1]
        \Input {LE point $x^* \in \mathcal X$, continuous MOP $F: \mathcal X \rightarrow \mathbb R^m$}
        \EndInput
        
        \Output{Approximation to LE set corresponding to $x^*$, \\
                \hspace{1.1cm} LE points belonging to superposed LE sets}
        \EndOutput
        \Procedure{ExploreEfficientSet}{$x^*, F$}

            \State $\texttt{set} \gets [x^*]$, $\texttt{superposed} \gets []$
            
            \For{$\text{obj} \in \{1, 2\}$}
                \State $\sigma \gets \sigma_\text{min}$
                \For{$t=1,2,\dots$}
                    \State Determine $x^{(t-1)}$, $x^{(t-2)}$, continuation direction $d$
                    \State $p \gets x^{(t-1)} + \sigma \cdot d / \|d\|$
                    \If{$f_\text{obj}(p) \geq f_\text{obj}^{(t-1)}$}
                        \If{$\sigma = \sigma_\text{min}$ \textbf{and} \texttt{use\_gradient}}
                            \State \textbf{break}
                        \EndIf
                        \State Update $\sigma, \texttt{use\_gradient}$
                        \State \textbf{continue}
                    \EndIf
                    \State $p^* \gets \Call{MultiObjectiveDescent}{p, F}$
                    \State $\phi \gets \Call{Angle}{x^{(t-2)} - x^{(t-1)}, x^{(t-1)} - p^*}$
                    \If{(($\|p - p^*\| > \sigma \textbf{ or } \phi > \phi_\text{max}$) \textbf{and} \\ \hspace{1.7cm} ($\sigma > \sigma_\text{min} \textbf{ or } \texttt{!use\_gradient}$))}
                        \State Update $\sigma, \texttt{use\_gradient}$
                        \State \textbf{continue}
                    \EndIf
                    \If{($(\|x^{(t-1)} - p^*\| > \sigma_\text{max} \textbf{ and } \texttt{use\_gradient})$ \textbf{ or } \\
                             \hspace{1.8cm} $ F(p^*) \prec F^{(t-1)}$)}
                        \State $\texttt{superposed}.append(p^*)$
                        \State \textbf{break}
                    \EndIf
                    \State Insert $p^*$ into \texttt{set}
                    \State Update $\sigma$, \texttt{use\_gradient}
                \EndFor
            \EndFor
            
            \State \textbf{return} $(\mathtt{set}, \mathtt{superposed})$
    \EndProcedure
    \end{algorithmic}
\end{algorithm}

\begin{table}[!t]
    \vspace*{-0.25cm}
    \caption{Suggested default values for the parameters of the set exploration algorithm (Alg.~\ref{alg:explore-set}). \texttt{diag} is defined as in Tab.~\ref{tab:parameters-descent}.}
    \vspace*{-0.2cm}
    \centering
    \begin{tabular}{|l|l||r||l|}
        \hline
        Parameter & Value & Description \\
        \hline
        \hline
        $\sigma_\text{min}$ & \texttt{1e-4} & Min. set step size\\
        $\sigma_\text{max}$ & \texttt{diag/100} & Max. set step size\\
        $\phi_\text{max}$ & \texttt{45}° & Max. angle per step\\
        $\lambda_\text{explore}$ & \texttt{2} & Set scaling factor\\
        \hline
    \end{tabular}
    \vspace*{-0.25cm}    
    \label{tab:parameters-explore-set}
\end{table}

From a (close to) LE starting point $x^*$, we explore the LE set into the direction of the first and second objective successively. Starting with the minimal step size %
$\sigma_\text{min}$, the PC method is executed as described above. For the first point in a set, %
the gradient of the currently traced objective will be used for the prediction, otherwise the secant through the two most recent points along the set %
is used.

After the prediction is performed with step size $\sigma$ along the predicted set direction, %
we ensure that we have improved the currently considered objective. 
If that is not the case, $\sigma$ is reduced by a factor of $\lambda_\text{explore}$, or, if $\sigma = \sigma_\text{min}$, we enforce the usage of the gradient direction and the iteration is restarted with the updated parameters. If this check fails while using the gradient direction, we terminate. %

Then, starting in the predicted solution $p$, a MO descent towards a LE point is performed. 
If the resulting descent traveled a longer distance than the original prediction (not using the gradient direction), the step size is reduced as above and the next iteration starts. %
In practice, we may terminate the MO descent early, if any of the accepted points along the search are more than $\sigma$ away from the starting point. %
Similarly, MOLE terminates, 
if the angle between the directions of the last step and the currently proposed step are larger than $\phi_\text{max}$, indicating that a smaller step size should be used.

Next, MOLE decides whether the obtained LE point belongs to the current set, or if MOLE crossed into another attraction basin. The latter is assumed, if either (1) a dominating point was found, or (2) the distance to the previous point on the set is greater than $\sigma_\text{max}$. %
Otherwise, we insert the new point into the current set and update $\sigma$ and $\texttt{use\_gradient}$, respectively.

At last, MOLE needs to decide whether a new LE point belongs to a previously explored set. In case MOLE finds such a set, the current point will be inserted into it and MOLE terminates. %
It recognizes a point as lying between two consecutive points $x^{(1)}$ and $x^{(2)}$ of a set, if the following conditions hold (\textsc{FindContainingSet}): In the objective space, it is located between the ideal and nadir points of $x^{(1)}$ and $x^{(2)}$, and in the decision space it is either closer to both $x^{(1)}$ and $x^{(2)}$ than those points are to each other, or closer to any of them than $\sigma_\text{min}$ (in this case, disregarding the objective space).

\subsection{Post-Processing}

In principle, a piece-wise linear model of the LE sets found by MOLE could be used as its output for a given application. Yet, to demonstrate its potential on conventional MO benchmarks, we incorporate a post-processing method that refines the solution by sampling further points in sparsely covered, though potentially highly relevant, areas along the LE sets. 
Here, we will use a procedure that is designed to optimize the dominated HV~\cite{zitzler2003performance} covered by all evaluated solutions, as this is the main target of the Bi-Objective BBOB \citep{tusar2016bbobbiobj} used in the upcoming experimental part. The overall idea is to first find all pairs of consecutive points in the LE sets whose ideal point is nondominated, i.e., they are likely to contain nondominated points between them. From all these pairs, at each iteration, the one with the largest HV gap between them is selected for sampling a new point at their mean in the decision space.

To ensure that the new solution is LE, we usually apply a MO descent. However, most LE sets are, after sufficient post-processing, modeled very well in the decision space using the piece-wise linear approximation inherent to MOLE.
\begin{figure}[!t]
    \centering
    \includegraphics[width=0.8\columnwidth]{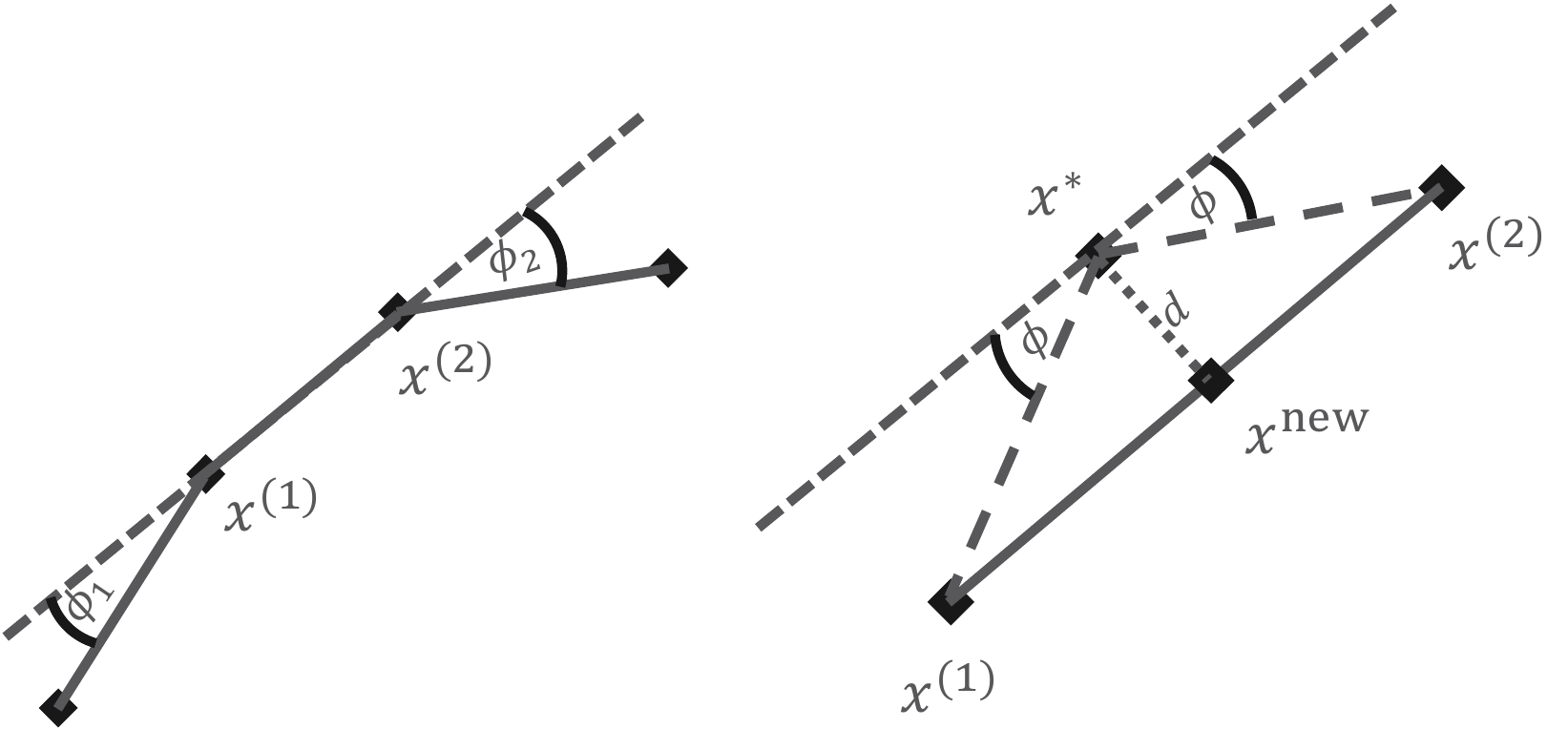}
    \vspace*{-0.2cm}
    \caption{
    Illustration of the \emph{maximal expected descent} $d$ for a new point between $x^{(1)}$ and $x^{(2)}$ with $\phi = \max(\phi_1, \phi_2)$.
    }
    \vspace*{-0.3cm}
    \label{fig:max-descent}
\end{figure}
This %
can be exploited as illustrated in Fig.~\ref{fig:max-descent}. 
Assuming that the set's model does not become worse during refining, 
i.e., the angles between consecutive steps in the set do not increase after inserting a new point, we can estimate the maximal distance $d$, which the MO descent would travel, by
\begin{equation}
    d = 0.5 \cdot \|x^{(1)} - x^{(2)}\| / \tan((180^\circ - \phi) / 2)
\end{equation}
with $\phi = \max\{\phi_1, \phi_2\}$. 
If $d$ is less than the minimal step size $\alpha_\text{min}$ of the MO descent, the descent is skipped, saving at least the function evaluations %
required for the initial gradient computation.
The resulting post-processing procedure is shown in Alg.~\ref{alg:post_process}. This process is repeated until a normalized HV target $0 < \theta \leq 1$ is reached.

\begin{algorithm}[!t]
\small
    \caption{Post-Processing LE Sets for HV Maximization}\label{alg:post_process}
    \begin{algorithmic}[1]
        \Input {LE sets \texttt{sets}, norm. HV target $\theta$, MOP $F: \mathcal X \rightarrow \mathbb R^m$}
        \EndInput
        
        \Output{Refined \texttt{sets} fulfilling (approx.) norm. HV target $\theta$}
        \EndOutput
        \Procedure{PostProcessHV}{$\texttt{sets}, \theta, F$}
            
            \State $\texttt{hv\_gap} \gets \text{Sum of HV gaps for } \texttt{sets}$
            
            \State $\texttt{nondominated} \gets \Call{GetNondominated}{\texttt{sets}}$
            
            \State $\texttt{max\_hv} \gets \text{Area between ideal, nadir of }\texttt{nondominated}$
            
            \While{$\texttt{hv\_gap} / \texttt{max\_hv} > \theta$}
                \State Find pair $x^{(1)}, x^{(2)}$ with maximal HV gap
                \State $x^\text{new} \gets (x^{(1)} + x^{(2)}) / 2$
                \State \texttt{max\_descent} $\gets$ Max. expected descent for $x^{(1)}, x^{(2)}$
                \If {$\texttt{max\_descent} > \alpha_\text{min}$}
                    \State $x^\text{new} \gets \Call{MultiObjectiveDescent}{x^\text{new}, F}$
                \EndIf
                \If {$\Call{Ideal}{x^{(1)}, x^{(2)}} \prec x^\text{new} \prec \Call{Nadir}{x^{(1)}, x^{(2)}}$}
                    \State Insert $x^\text{new}$ into corresponding $\texttt{set}$
                \EndIf
                \State Update $\texttt{hv\_gap}$
            \EndWhile

            \State \textbf{return} $\texttt{sets}$
    \EndProcedure
    \end{algorithmic}
\end{algorithm}

\section{Experimental Study and Discussion}
\label{sec:exp_discuss}

After having introduced MOLE in detail, this section evaluates its potential. Keep in mind that MOLE remains a local optimizer that is not expected to be regularly competitive to global meta-heuristics. %

\begin{figure*}[!t]
    \centering
    \includegraphics[width=0.19\textwidth]{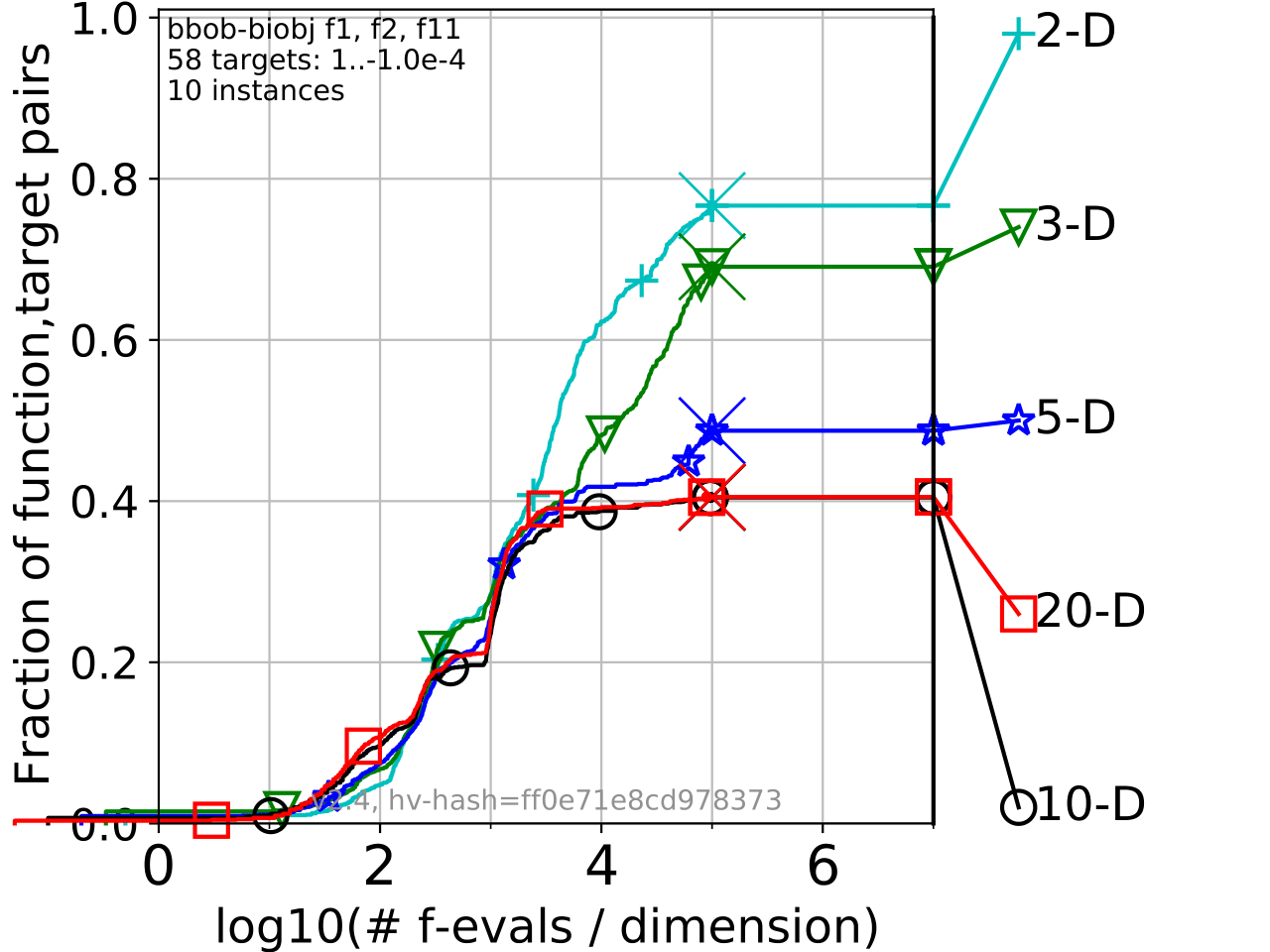}
    \includegraphics[width=0.19\textwidth]{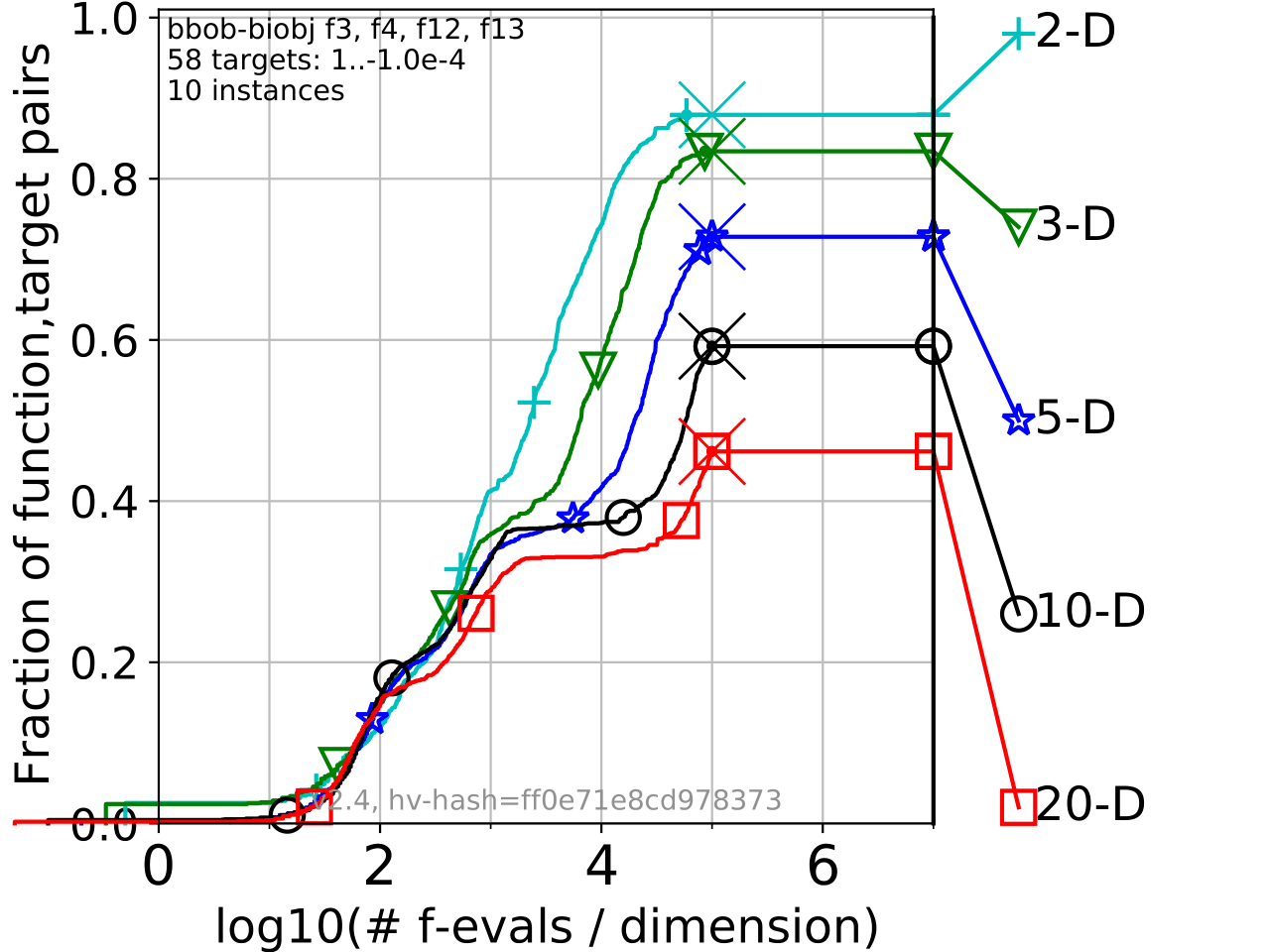}
    \includegraphics[width=0.19\textwidth]{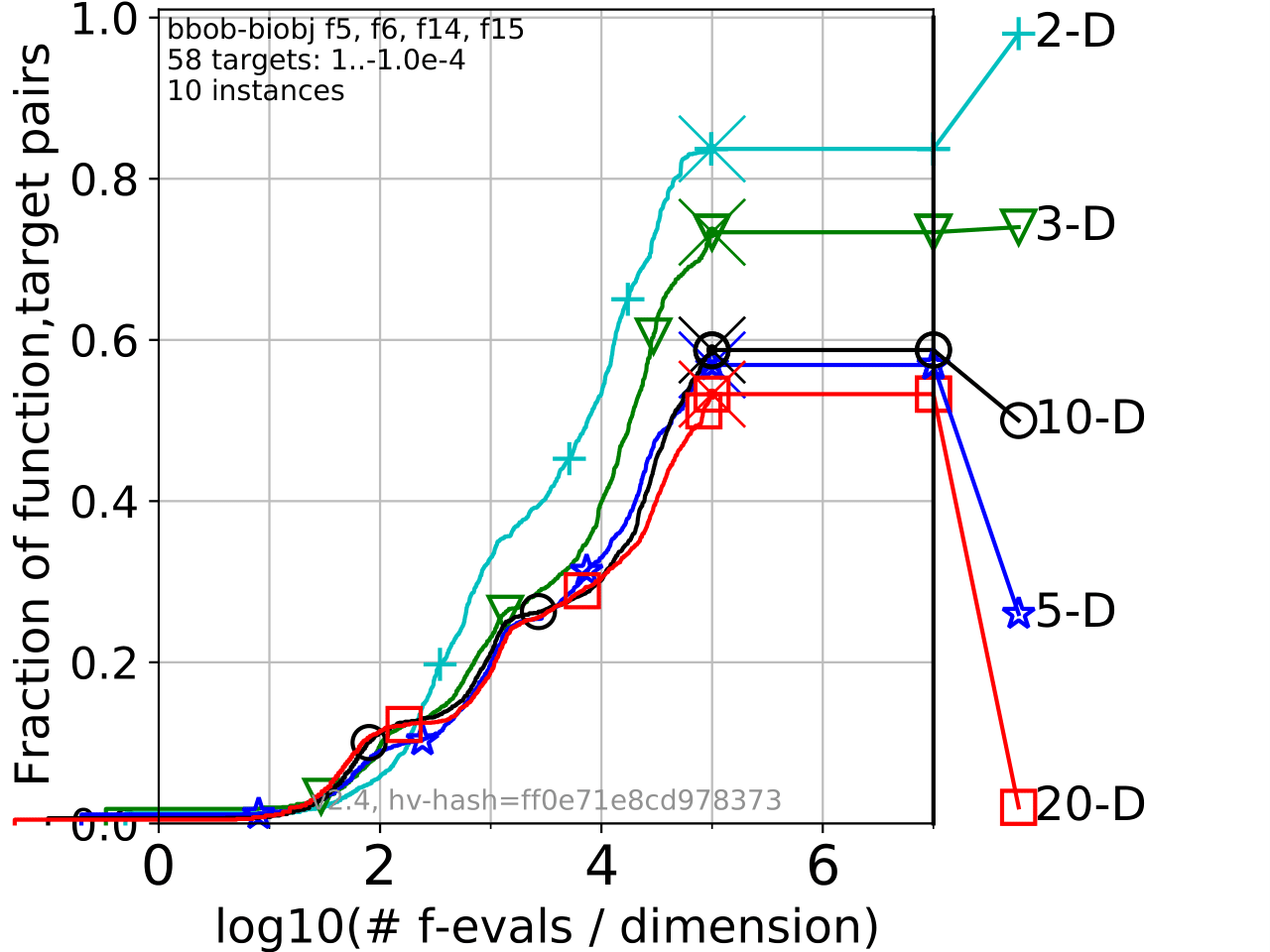}
    \includegraphics[width=0.19\textwidth]{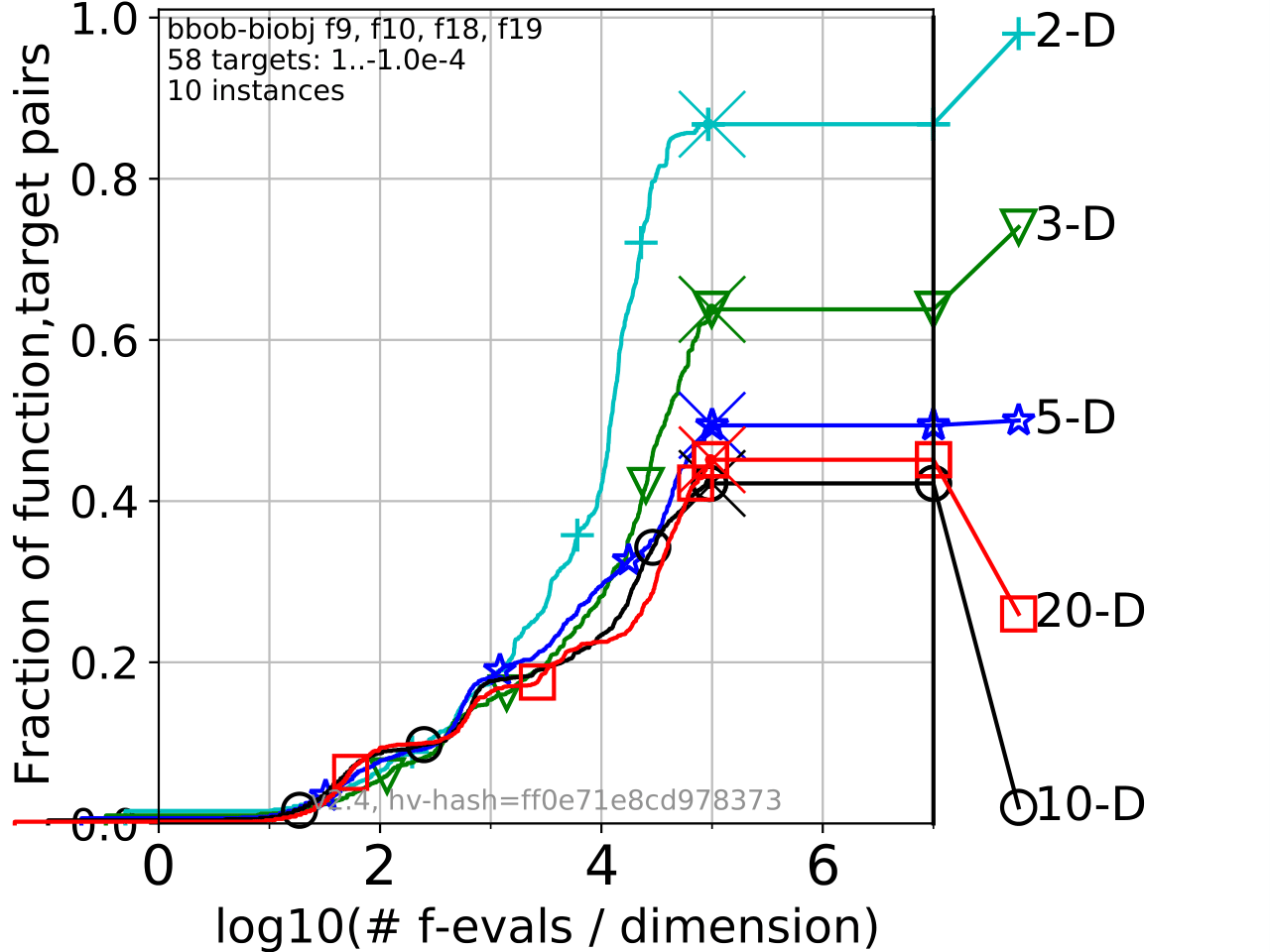}
    \includegraphics[width=0.19\textwidth]{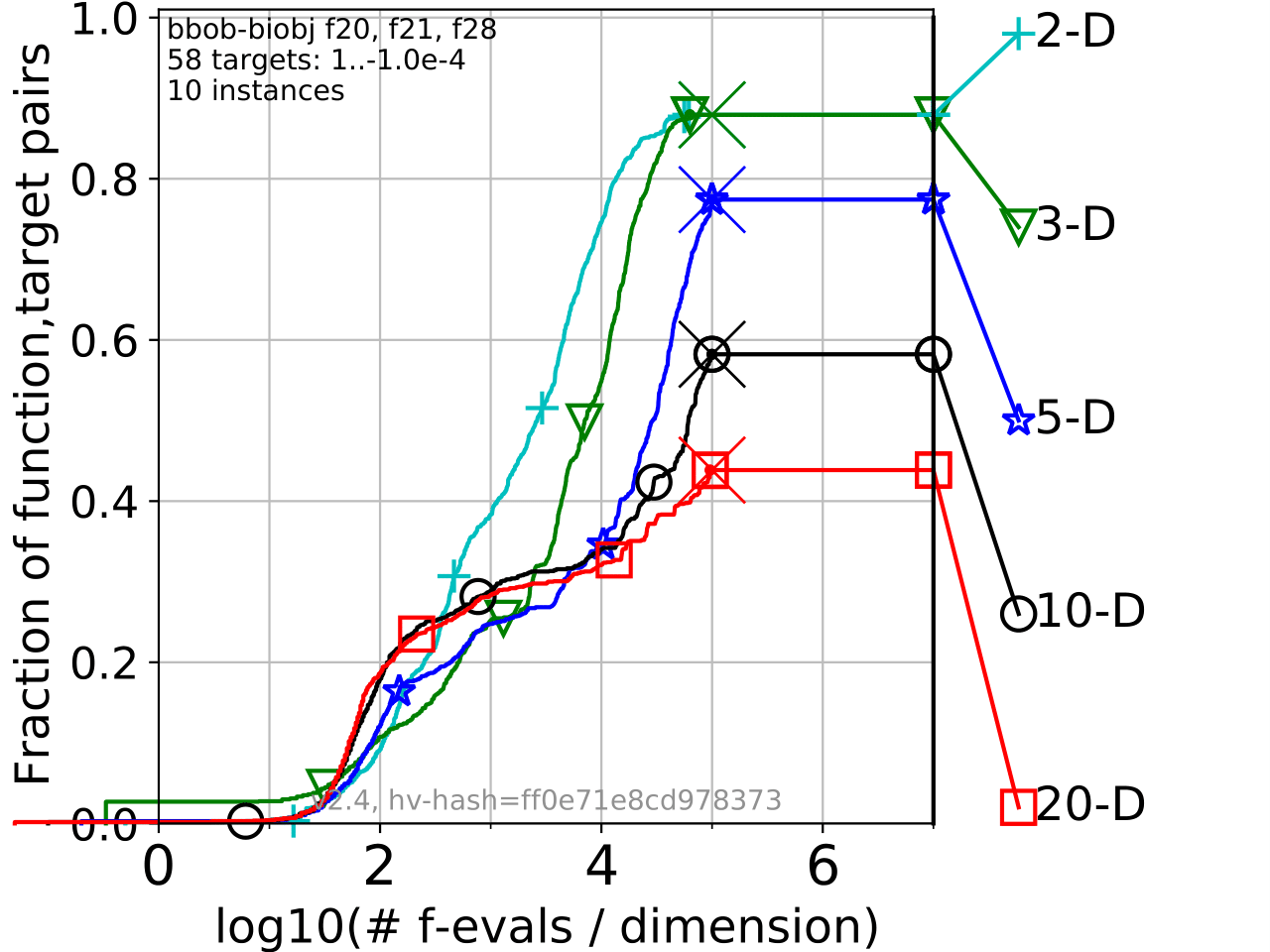}
    \includegraphics[width=0.19\textwidth]{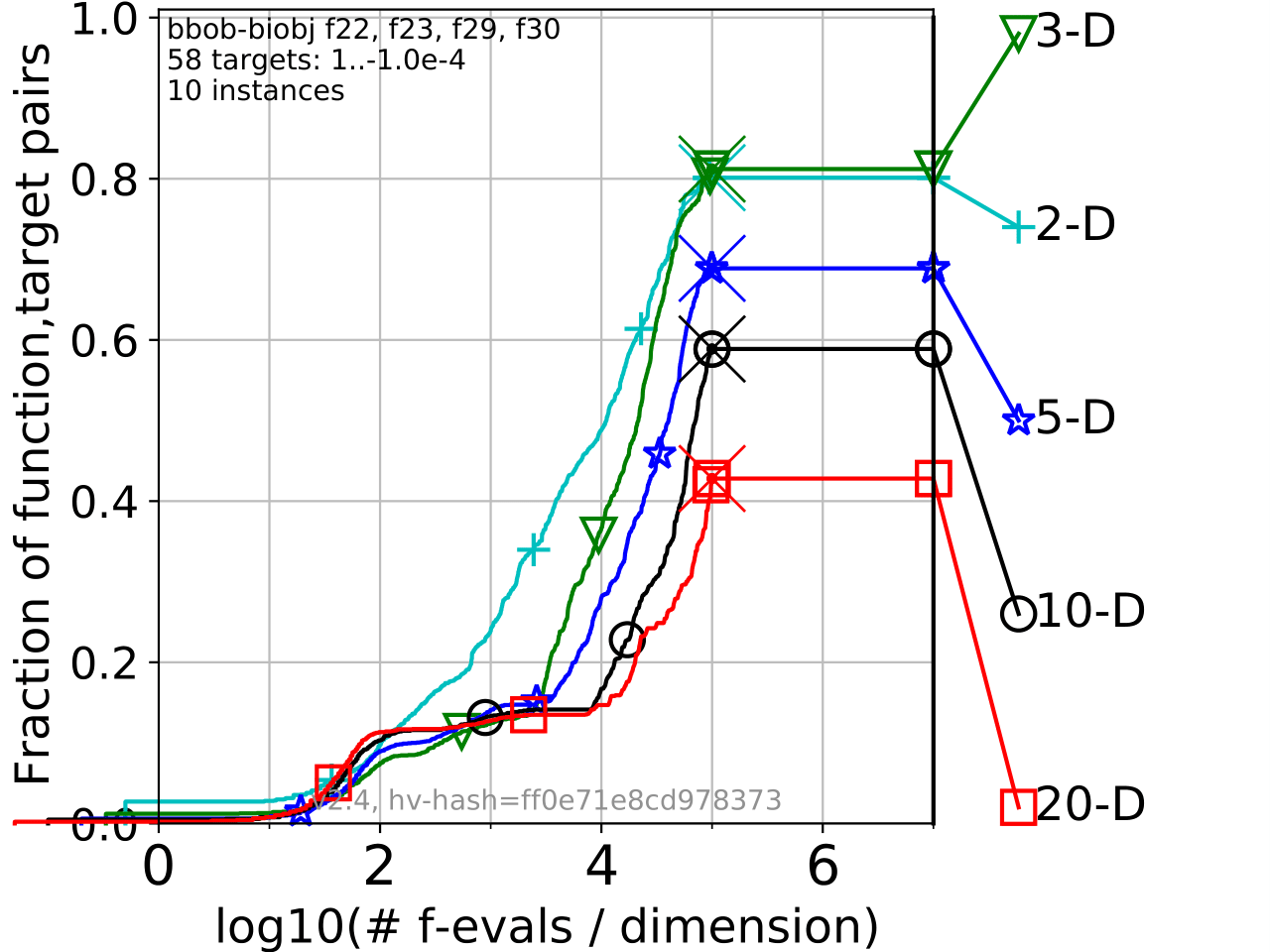}
    \includegraphics[width=0.19\textwidth]{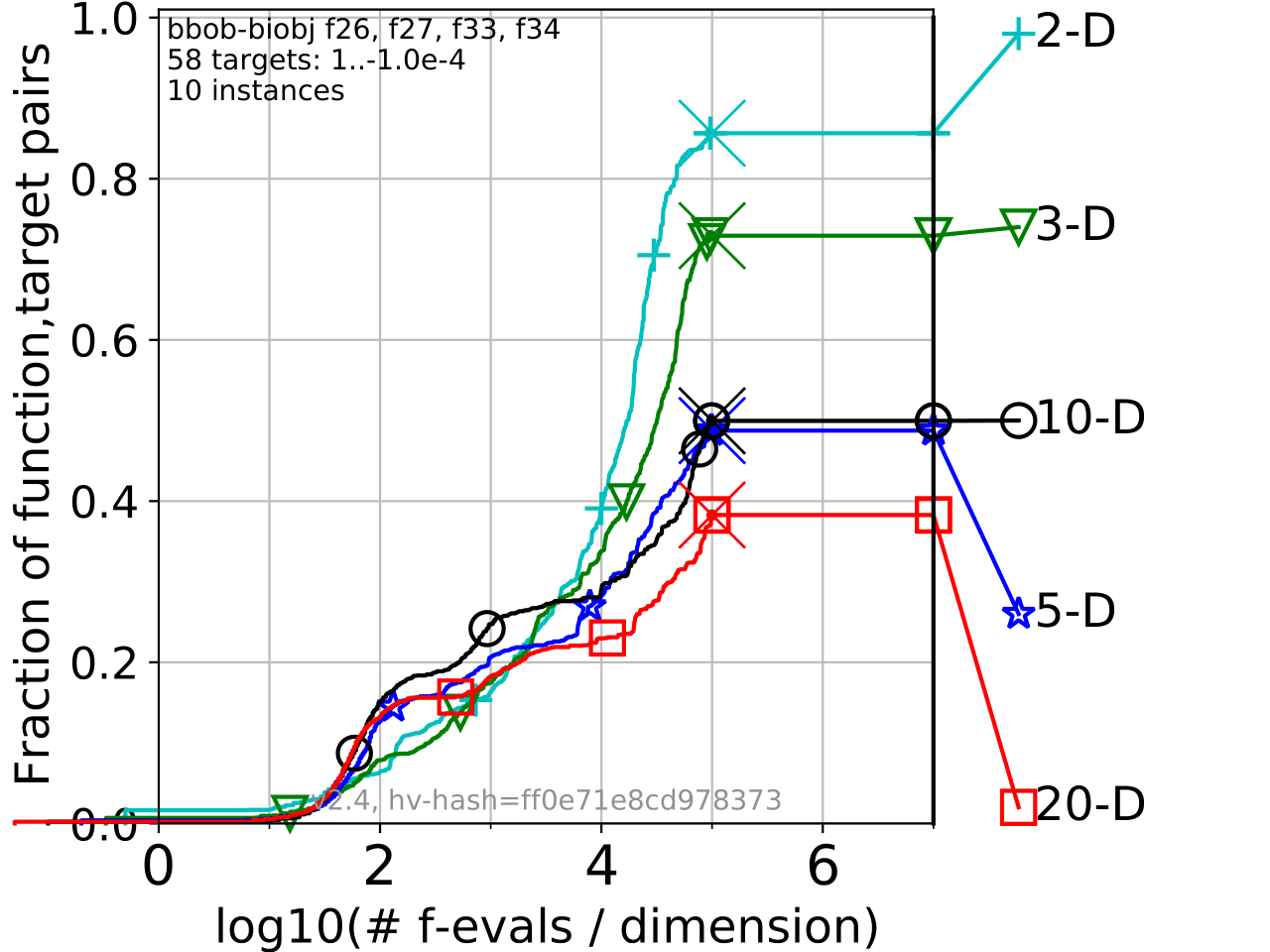}
    \includegraphics[width=0.19\textwidth]{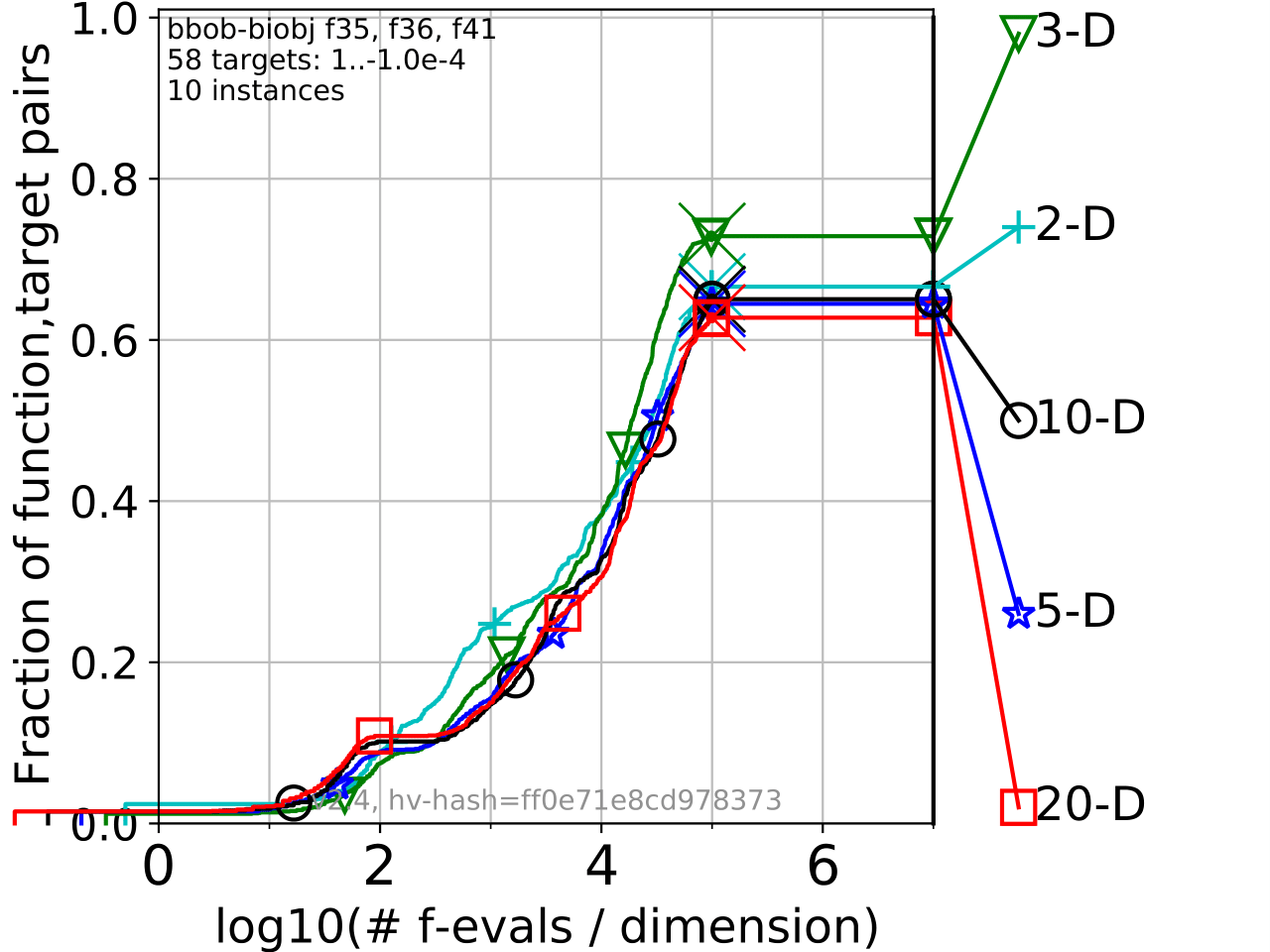}
    \includegraphics[width=0.19\textwidth]{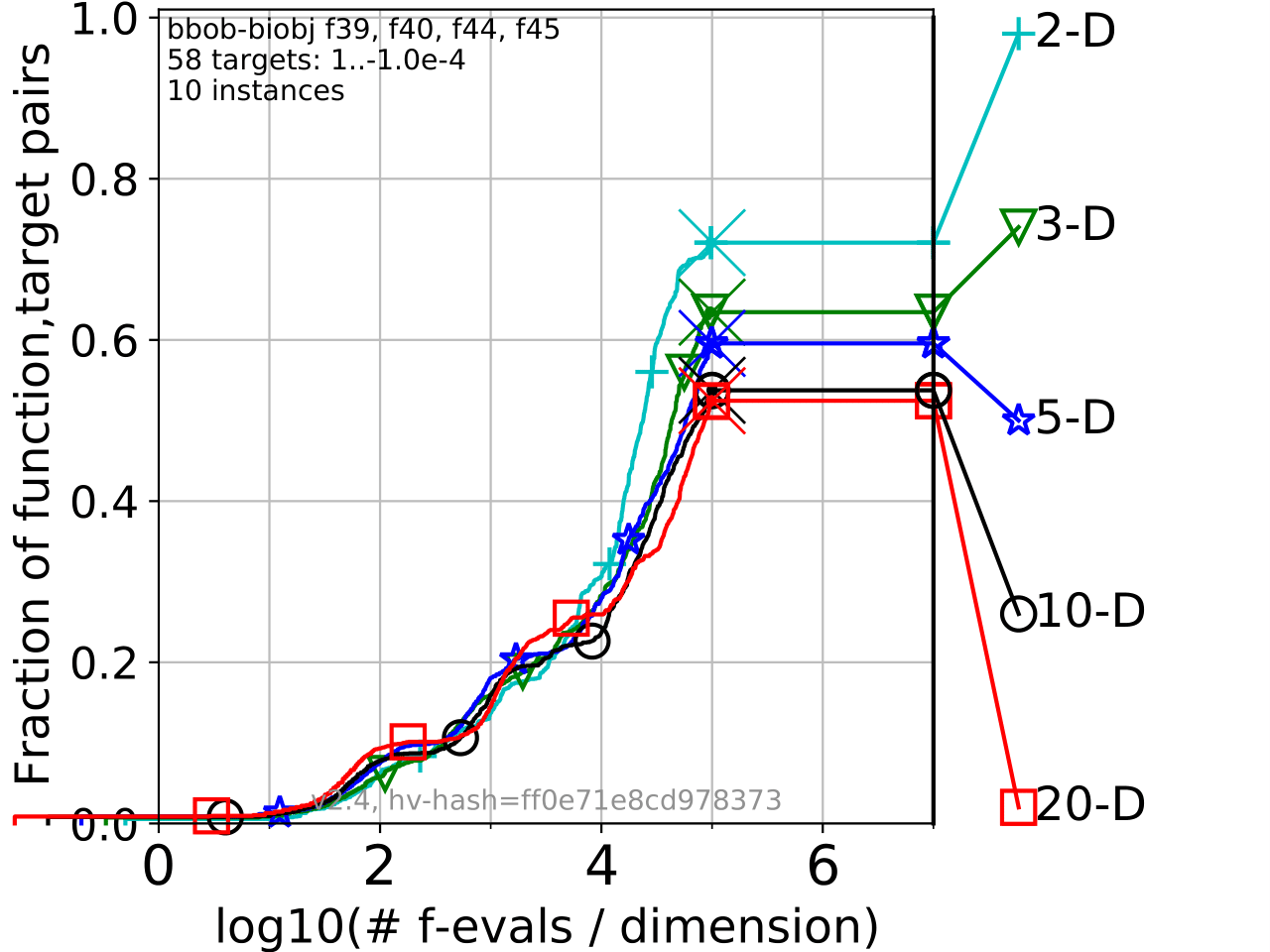}
    \includegraphics[width=0.19\textwidth]{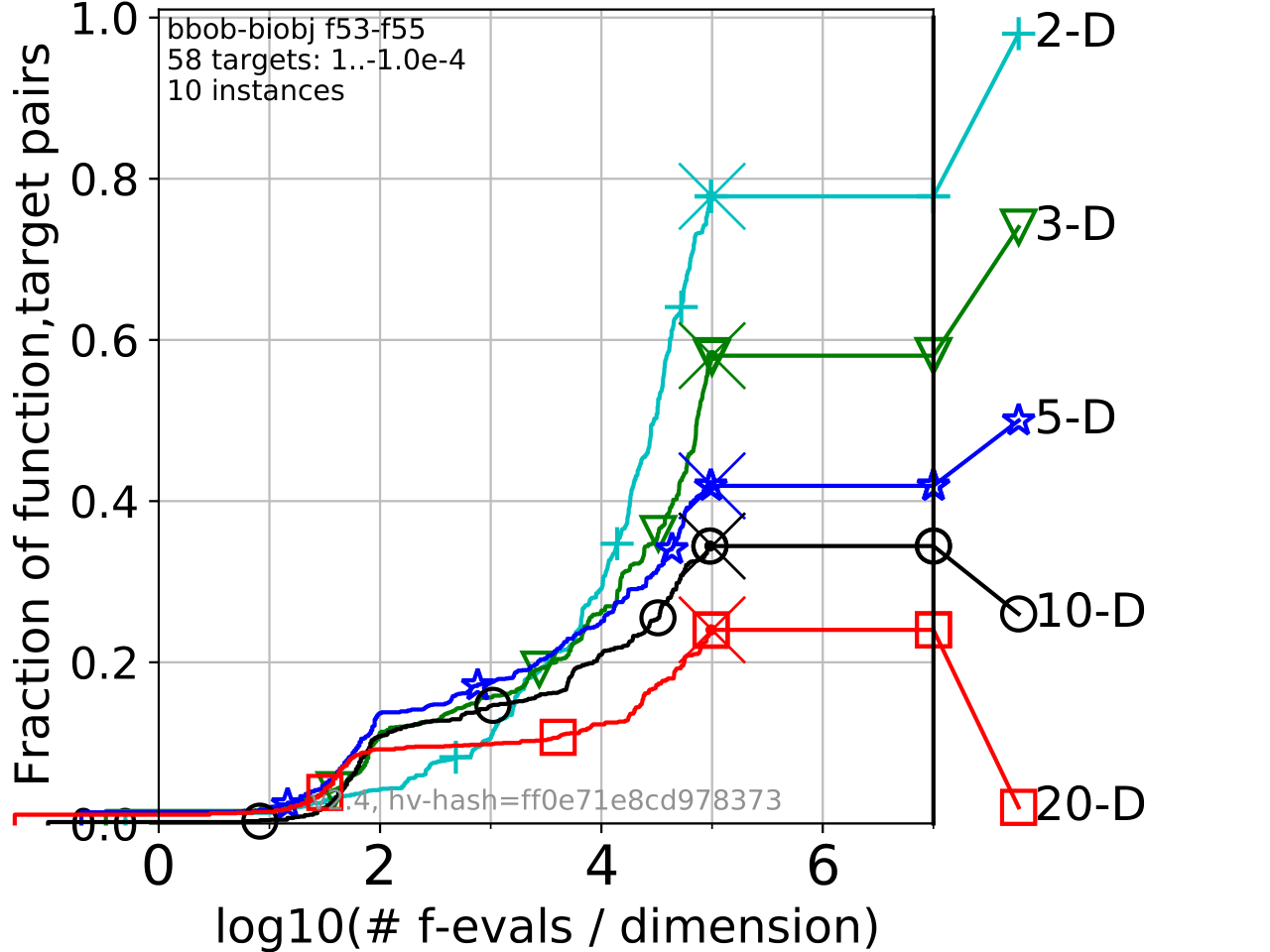}

    \includegraphics[width=0.19\textwidth]{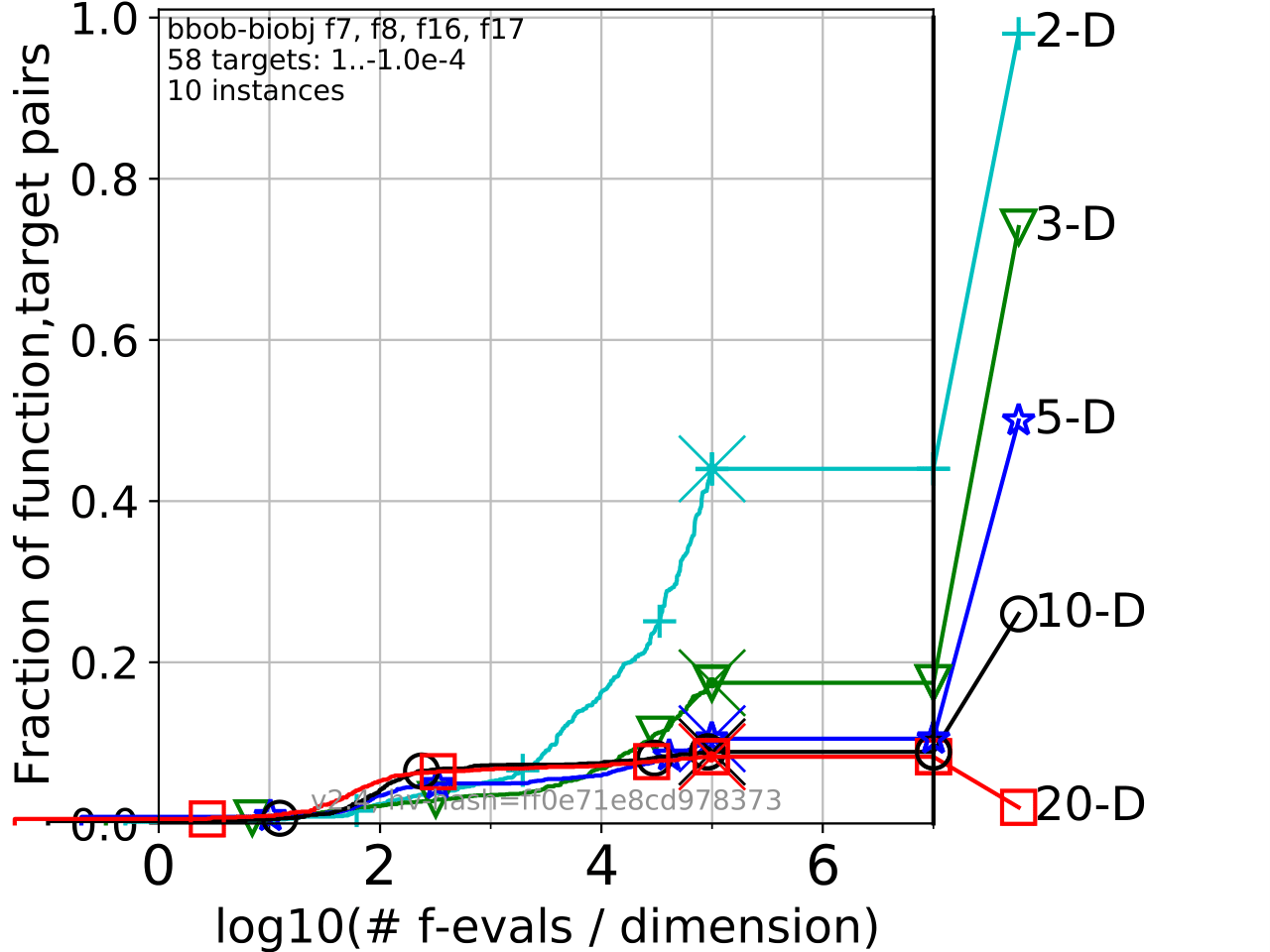}
    \includegraphics[width=0.19\textwidth]{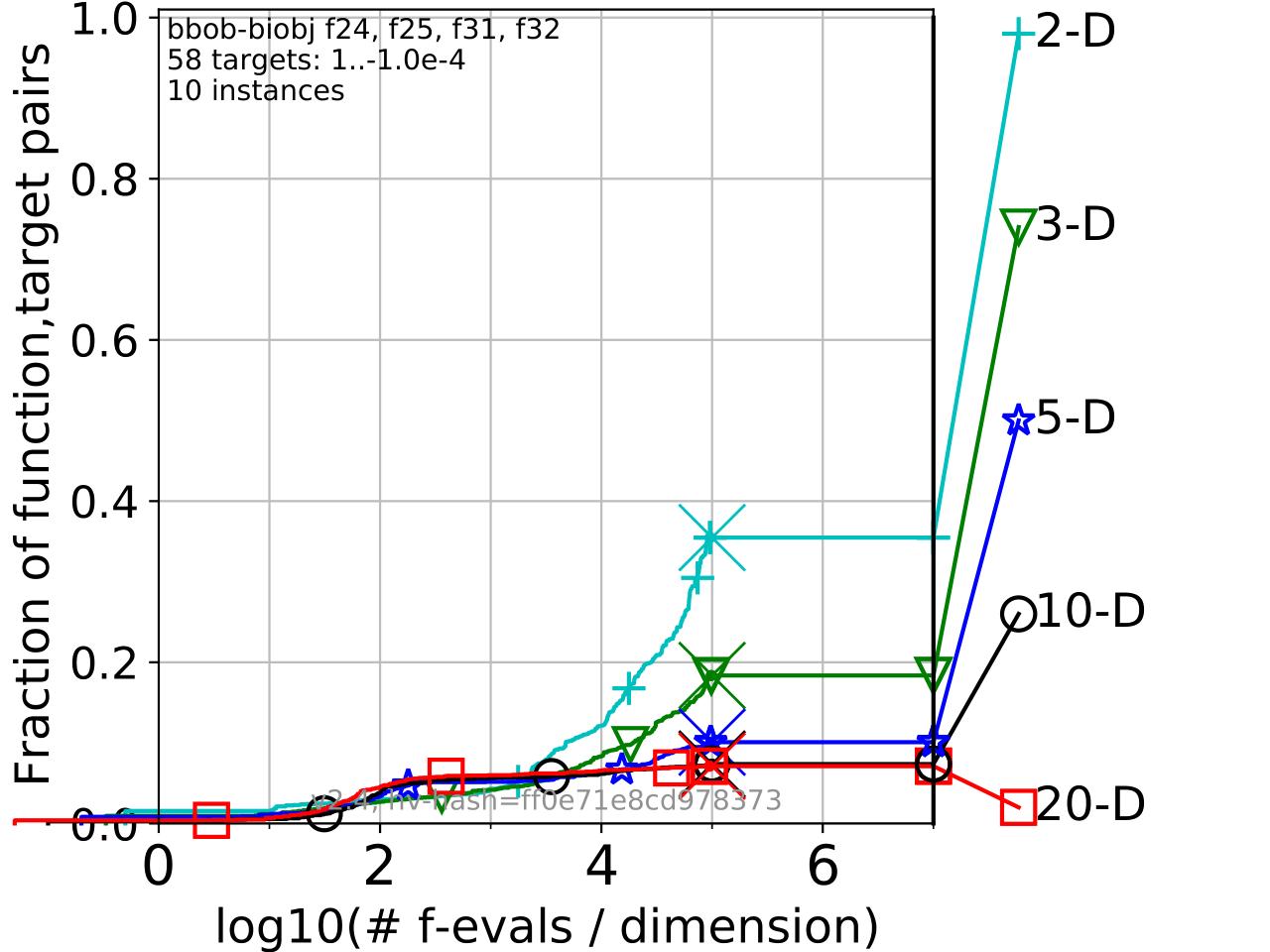}
    \includegraphics[width=0.19\textwidth]{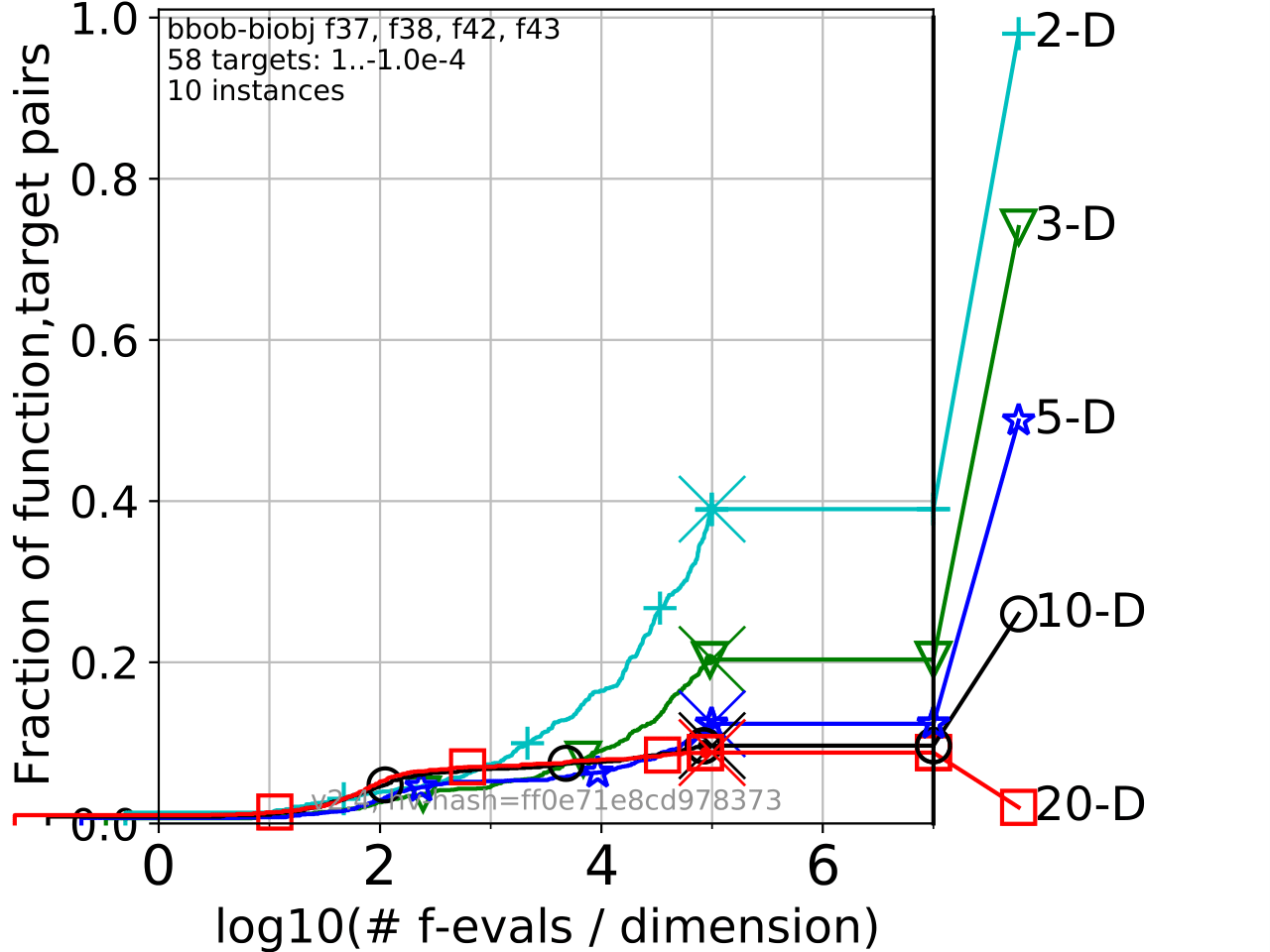}
    \includegraphics[width=0.19\textwidth]{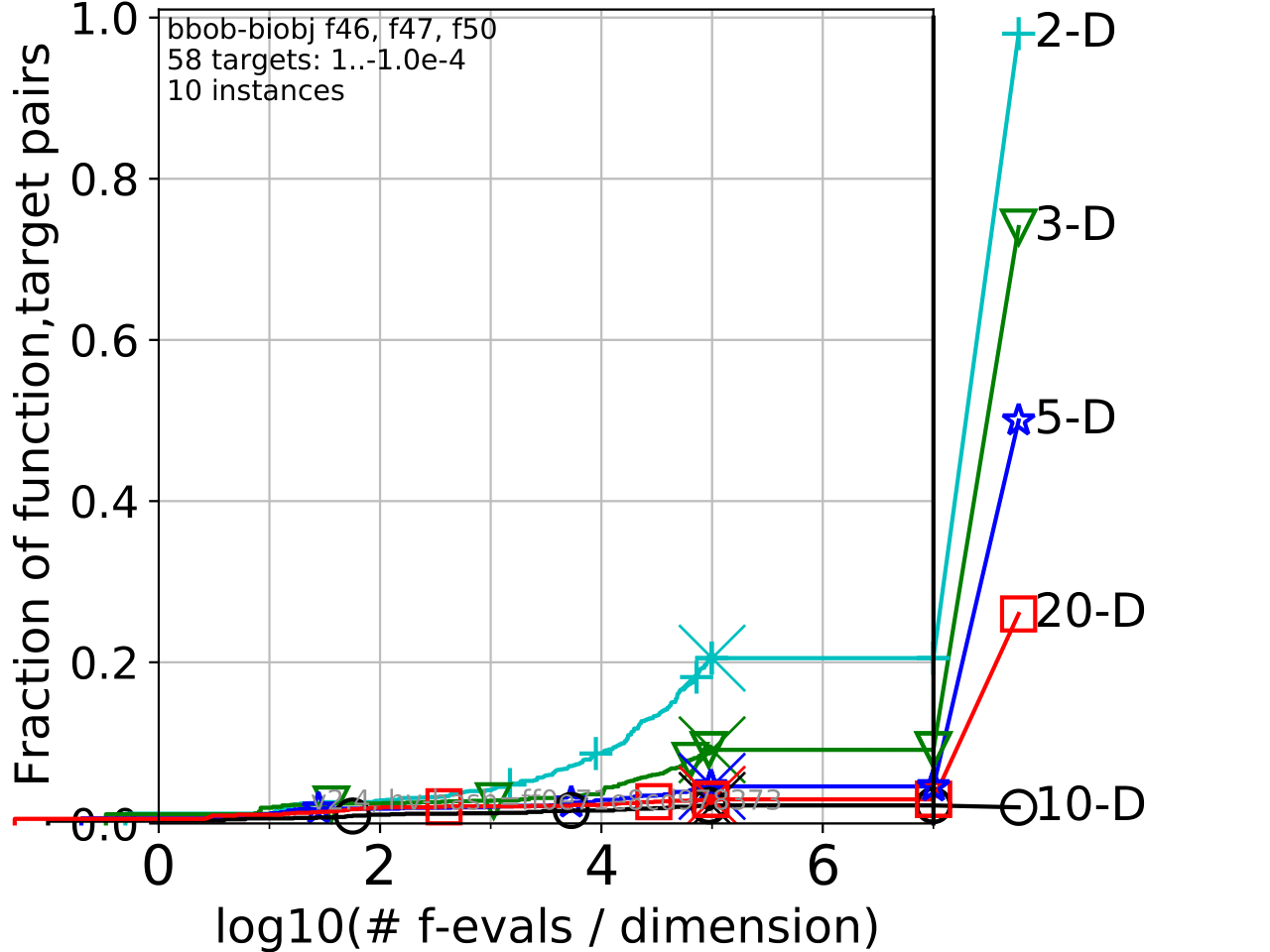}
    \includegraphics[width=0.19\textwidth]{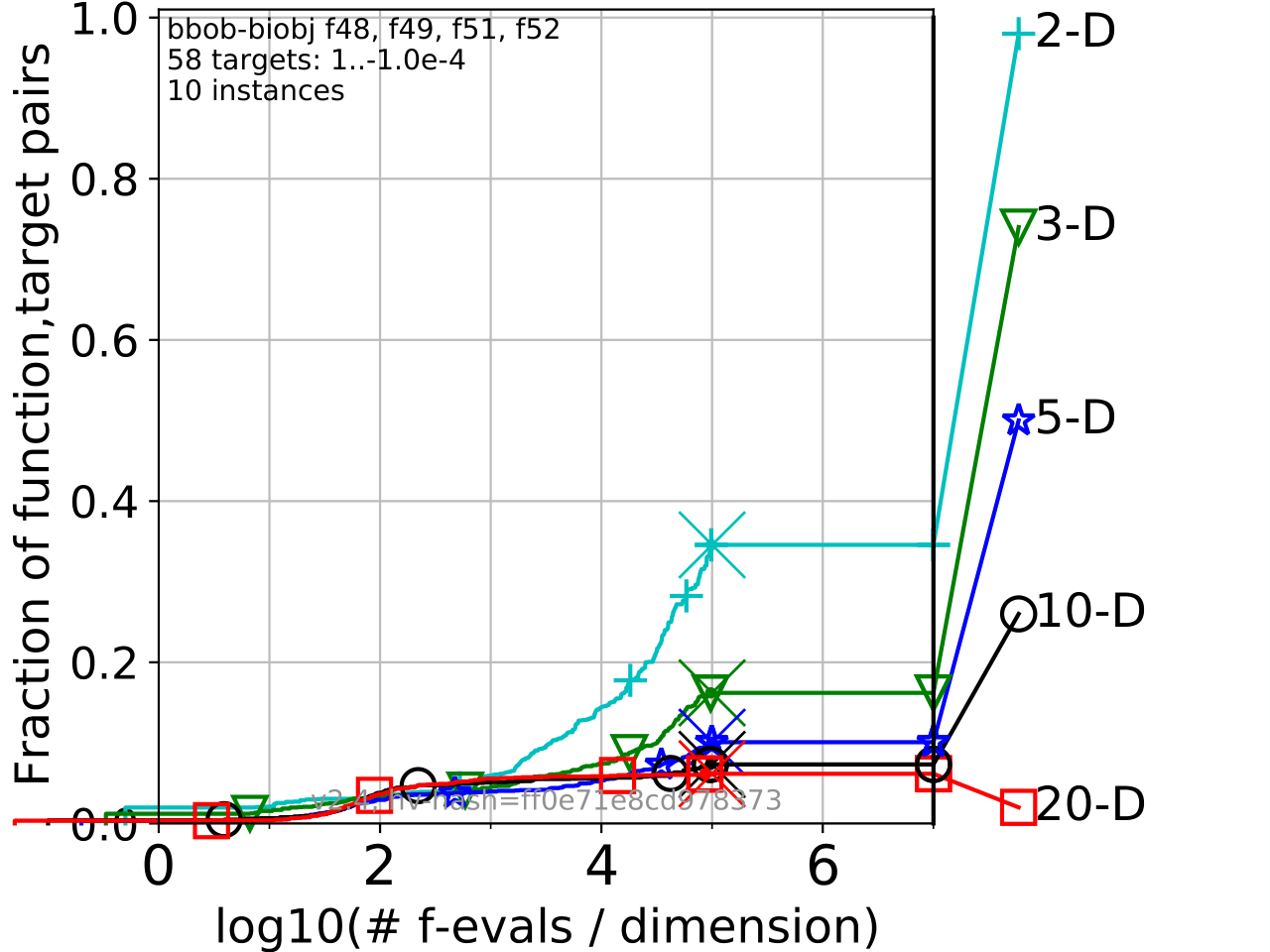}
    \vspace*{-0.25cm}
    \caption{
    ECDF graphs showing MOLE's performance (dimensions $2, 3, 5, 10$ and $20$) across groups of bi-objective BBOB functions.
    \vspace*{-0.2cm}
    }    
    \label{fig:coco-mole-groups}
    \centering
    \includegraphics[width=0.19\textwidth]{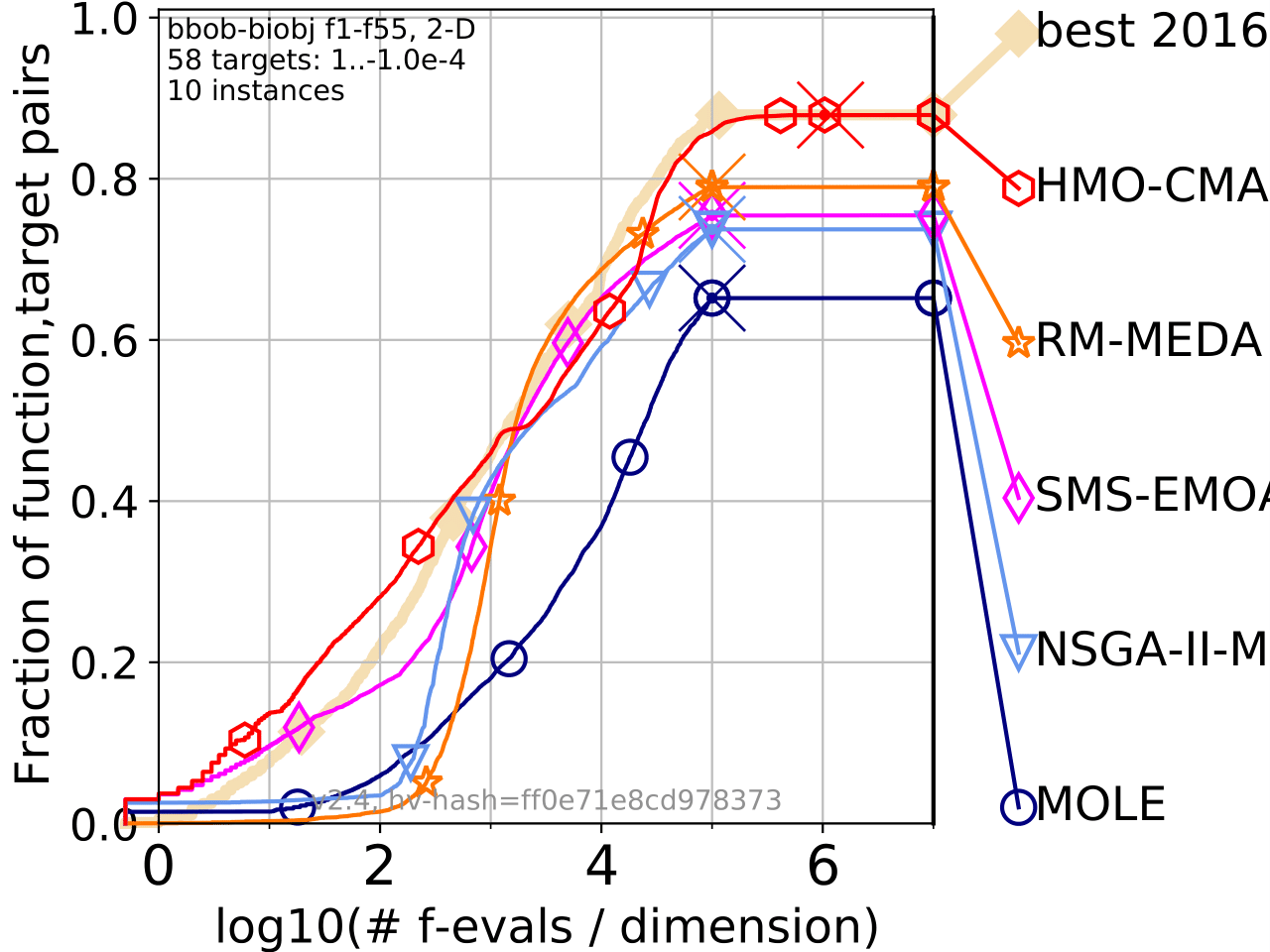}
    \includegraphics[width=0.19\textwidth]{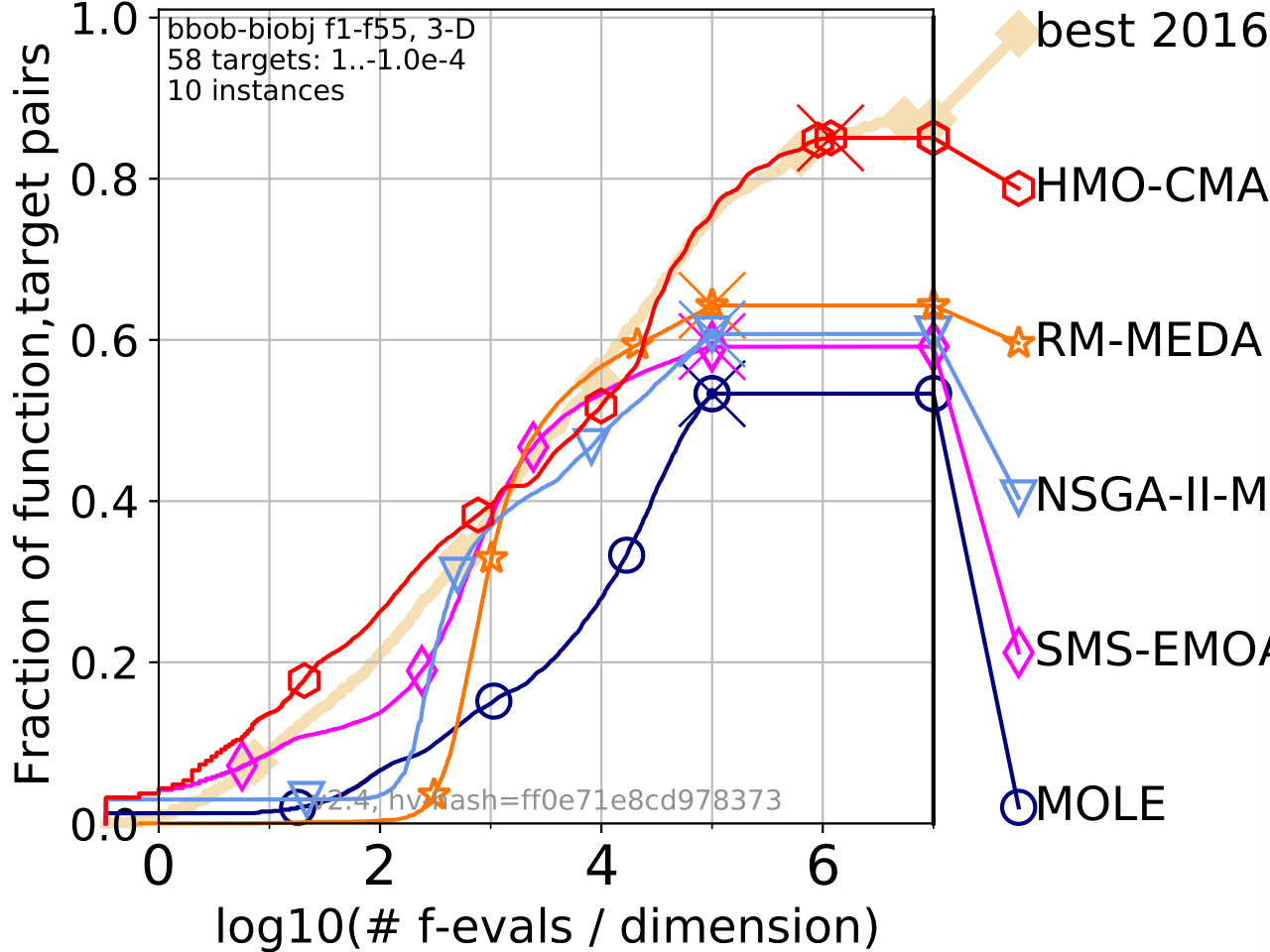}
    \includegraphics[width=0.19\textwidth]{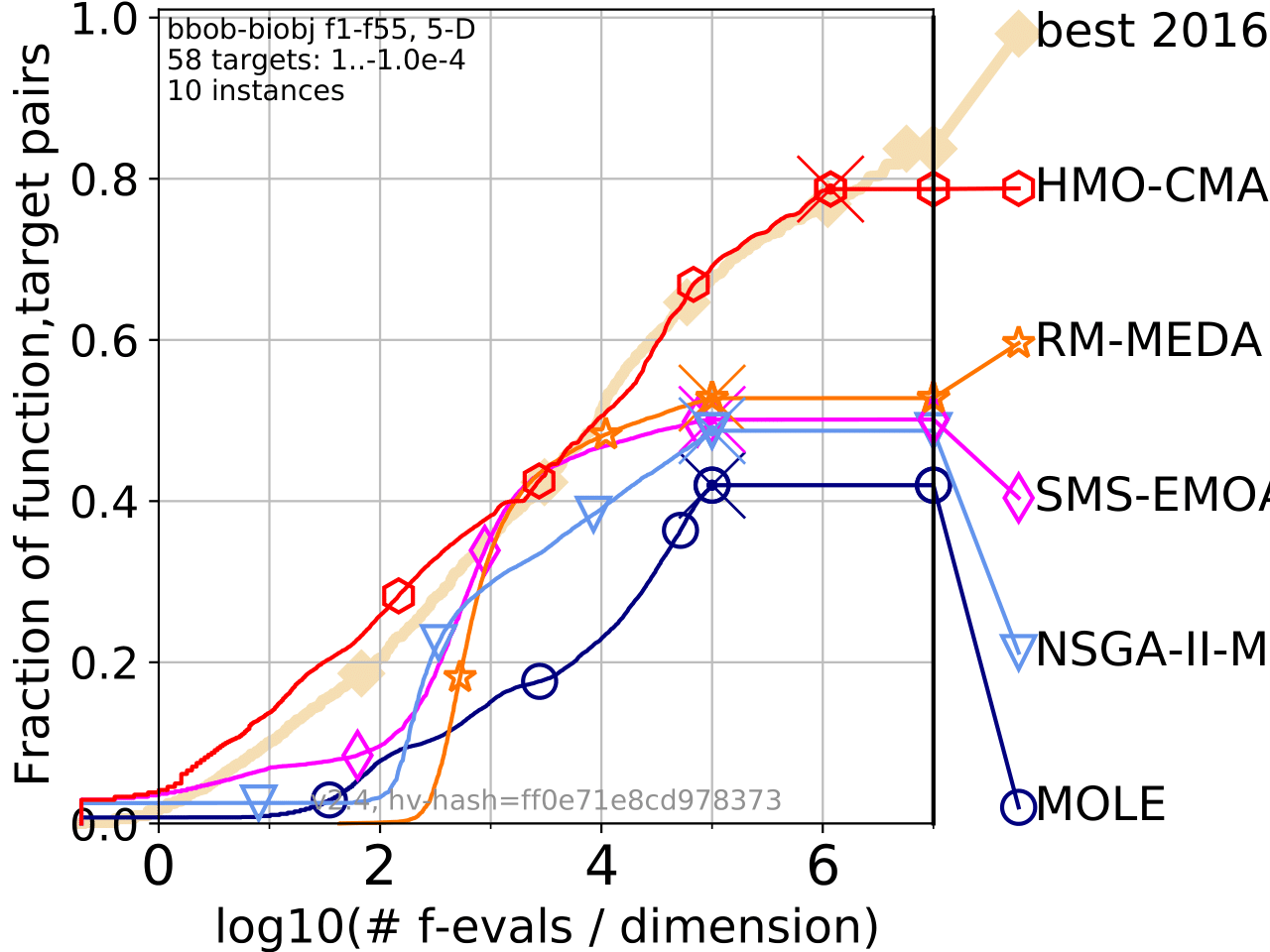}
    \includegraphics[width=0.19\textwidth]{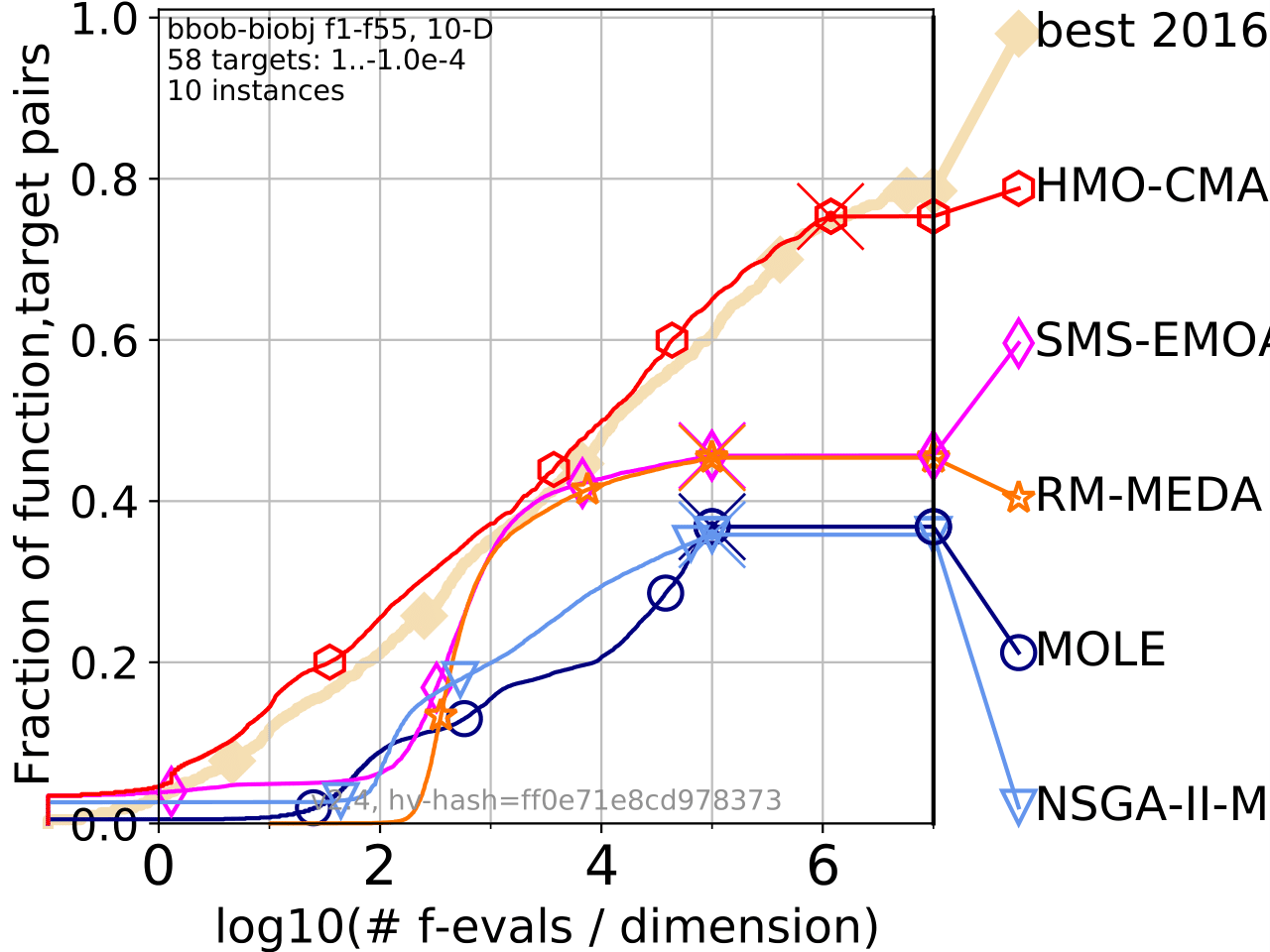}
    \includegraphics[width=0.19\textwidth]{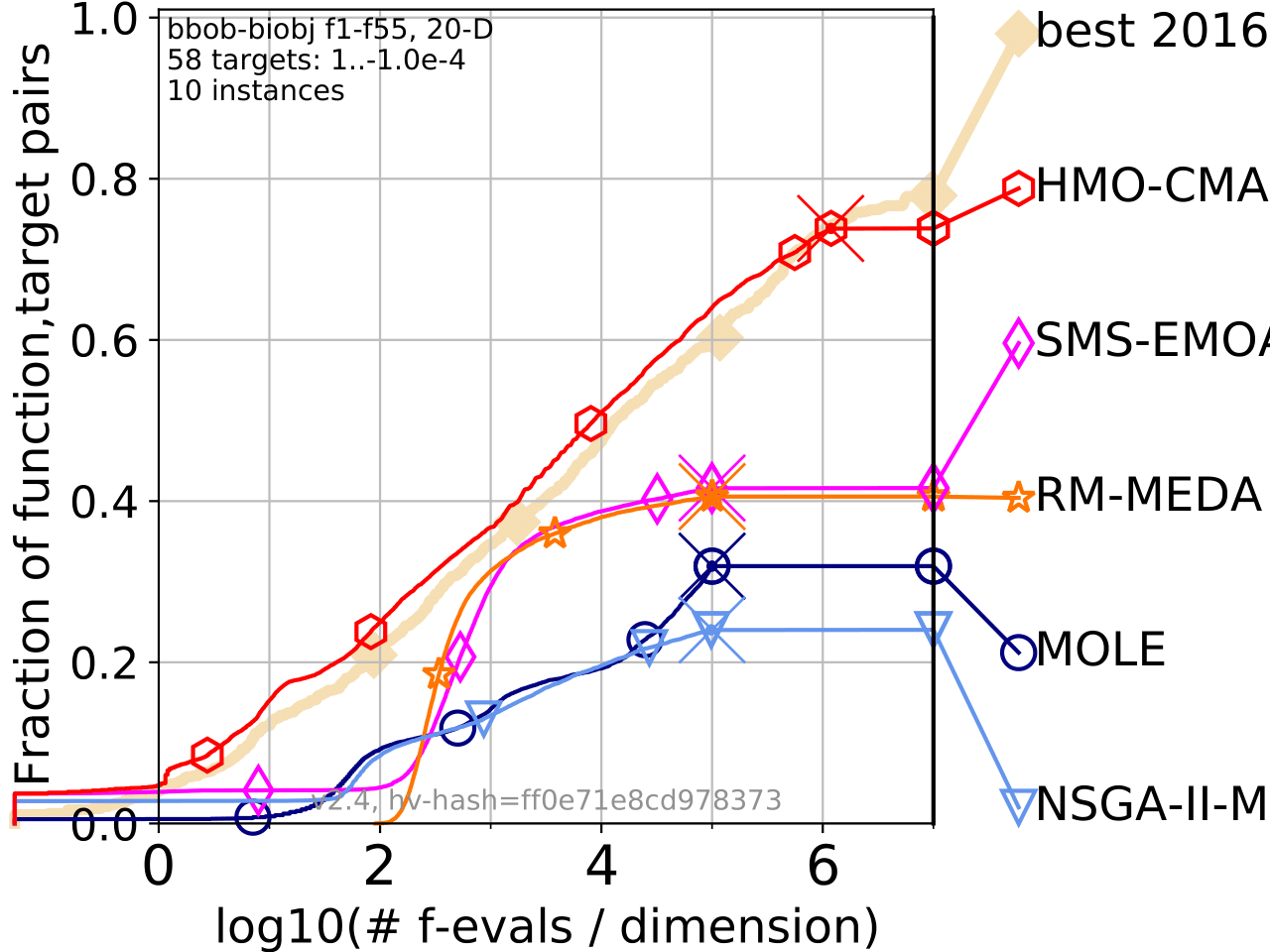}
    \vspace*{-0.25cm}
    \caption{ECDF graphs comparing MOLE's performance with several other state-of-the-art optimizers, and the Best-2016 reference algorithm on all $55$ functions in dimensions $2,3,5,10$ and $20$.}
    \vspace*{-0.25cm}
    \label{fig:coco-comparison}
\end{figure*}

\subsection{Experimental Setup}
We benchmark our bi-objective MOLE using the Bi-Objective BBOB \citep{tusar2016bbobbiobj} as it presents one of the most prominent test suites. %
The COCO platform allows to integrate benchmark data of other algorithms and enables the comparison between different optimizers. This ensures that all optimizers are configured well for the benchmark, which should enable a fair comparison of the algorithms.

As quality indicator, the normalized dominated HV w.r.t. the ideal and nadir points (in objective space) is used if any evaluated point dominates the nadir point, and the (negative) distance to the dominating area, otherwise. 
To generate similar targets as for the SO BBOB functions, a reference HV value needs to be known per problem. Yet, the optimal normalized HV values are not known for most benchmark functions. %
Instead, they are approximated based on an archive of runs from several optimizers, so that an optimizer could find a solution set with a better dominated HV. %
For this reason, some negative target values are introduced in addition to the usual positive ones, in order to reward optimizers that can find better solution sets. This results in the targets $\{-10^{-4},-10^{-4.2},\dots,-10^{-5},0,10^{-5},10^{-4.9},\dots,10^{0}\}$ %
being used by default, where only $52/58 \approx 90\%$ can definitely be reached. Even the $0$ target might be unobtainable for lower budgets, as this can require hundreds of thousands of nondominated points to solve.

In addition to evaluating MOLE alone, %
its performance will be compared to several state-of-the-art solvers from COCO's data archive. The selected solvers include the HMO-CMA-ES~\cite{loshchilov2016anytime}, a hybrid CMA-ES algorithm that is the single best solver on the Bi-Objective BBOB, the best performing variant of SMS-EMOA~\cite{BNE2007smsemoa}, the original NSGA-II~\cite{Deb2002nsg2} and the RM-MEDA~\cite{zhang2008rm} algorithm that also attempts to model LE sets, albeit embedded in an evolutionary technique. Finally, the Best-2016 reference algorithm is used to compare against the best performance per problem from the 2016 competition. %

With the Bi-Objective BBOB in mind, we can describe the configuration chosen for benchmarking MOLE. We run MOLE with the default parameters and a budget of $10^5 d$ on all $55$ functions and dimensions $2,3,5,10$ and $20$. The instances $1-10$ are chosen, as there is only reference data for these.
Further, we run MOLE for at most $1000$ starting points, which are sampled uniformly at random from $[-5,5]^d$, and terminate the search for further LE sets after $1000$ sets are discovered to limit the overhead on very multimodal functions.
The HV maximization post-processing is first used after $10$ starting points were successfully evaluated, in order to explore the landscape to some extend before optimizing the globally efficient regions found so far. It is then started each time after a new LE set that contains nondominated solutions is found. Finally, we set the target precision to $\theta = 2 \cdot 10^{-5}$, which resembles the final positive target value of the benchmark the closest. %

\subsection{Experimental Results}

The ECDF graphs in Fig.~\ref{fig:coco-mole-groups} summarize MOLE's performance across numerous groups of (similar) Bi-Objective BBOB functions. %
The source code for MOLE and an accompanying R-package are published at \textcolor{blue}{\url{https://github.com/schaepermeier/moleopt}}.
The full benchmarking data is available from \textcolor{blue}{\url{https://github.com/schaepermeier/2022-gecco-mole-data}}.
In general, the performance of MOLE worsens with increasing decision space dimensionality. Taking a closer look at the results of the individual function groups, however, reveals that much of the performance is lost on the multimodal functions with adequate global structure (bottom row in Fig.~\ref{fig:coco-mole-groups}). It turns out that these problems have properties adverse to MOLE, resp.~local search methods in general, which is discussed in more detail below.
While the performance does not decrease nearly as much on the other functions (top two rows in Fig.~\ref{fig:coco-mole-groups}), some performance loss can also be observed. %

Fig.~\ref{fig:coco-comparison} compares MOLE with several other state-of-the-art MO optimizers (for each evaluated dimension). 
The other algorithms with a maximum budget of $10^5 d$ seem to have converged around their budget cap, reaching a plateau of solved HV targets. In contrast, MOLE, as already seen above, did not fully converge on all problems at this point and would likely continue to solve HV targets if given further evaluations. Only the hybrid HMO-CMA-ES outperforms MOLE on most benchmark functions.

\subsection{Discussion}

The benchmarking results indicate that MOLE, as a pure local search algorithm, can compete on many problems with the established evolutionary algorithms that are prevalent in continuous MOO. In particular, functions that are moderately multimodal and moderately conditioned (or simpler) are efficiently optimized by it.

However, MOLE encounters problems when its core assumptions on the underlying landscape are violated. %
Either there are too many LE sets to efficiently track and explore, or weakly dominated areas complicate the modelling of the LE sets and the application of MO descent procedures. While the performance does not degenerate as much on the remainder of the test functions, increases in decision space dimensionality also lead to some performance losses.

\section{Conclusions}
\label{sec:conclude}

In this paper, we introduced the Multi-Objective Landscape Explorer (MOLE) as new local search approach. It is based on a prototype idea called MOGSA but exploits local properties of MO landscapes in a more sophisticated way. 
The potential for the exploitation of the local search structures was evaluated in a benchmarking study. We showed that MOLE is a competitive optimizer on a set of test problems in the Bi-Objective BBOB \citep{tusar2016bbobbiobj}. In particular, unimodal and weakly structured functions were optimized efficiently, regardless of dimension, competing with established algorithms like NSGA-II, even though all experiments were conducted with a rather simple configuration of MOLE. This includes functions that are multimodal in the sense of containing multiple, disconnected locally efficient sets, but are unimodal w.r.t.~the reachability of the globally efficient solution based on the set interactions. This does not only suggest MOLE as an alternative approach but also highlights the potential of this local search for integration in meta-heuristics. %

In future work, the default parameter values of MOLE should be optimized to evaluate their impact and improve the overall performance of the algorithm.
Further, hybridization of MOLE with an evolutionary algorithm (EA) might present a remedy for the problems MOLE encounters on certain highly multimodal functions.
The global search strengths of EAs could be combined with the exploitation of the local search dynamics in the globally efficient regions to refine the solutions proposed by the EA.

\subsection*{Acknowledgments}
The authors acknowledge support by the \href{https://www.ercis.org}{\textit{European Research Center for Information Systems (ERCIS)}}.

\bibliographystyle{unsrt}
\bibliography{arxiv}

\end{document}